\documentclass[english]{article}
\usepackage[utf8]{inputenc}
\usepackage{fullpage} 

\usepackage{amsmath}
\usepackage{amssymb}
\usepackage{graphicx}
\usepackage{subfigure}
\usepackage[dvipsnames]{xcolor}
\usepackage{epsfig}
\usepackage{epstopdf}
\usepackage{hyperref}
\usepackage{algorithm,algcompatible}
\usepackage{authblk}

\newcommand{\R}{\mathbb{R}}

\newcommand{\M}{\mathcal{M}}
\newcommand{\TM}{\mathcal{T}_f\mathcal{M}}
\newcommand{\ee}{\mathrm{e}}
\newcommand{\dd}{\, \mathrm{d}}

\graphicspath{ {./figures/} } 

\title{A low-rank projector-splitting integrator for the Vlasov--Maxwell equations with divergence correction}

\author[1]{Lukas Einkemmer}
\author[2]{Alexander Ostermann}
\author[3]{Chiara Piazzola}
\affil[1,2,3]{Department of Mathematics, University of Innsbruck, Austria \protect \\
\nolinkurl{lukas.einkemmer, alexander.ostermann, chiara.piazzola@uibk.ac.at}}

\date{\today}

\begin{document}

\maketitle

\section*{Abstract}
The Vlasov--Maxwell equations are used for the kinetic description of magnetized plasmas. As they are posed in an up to 3+3 dimensional phase space, solving this problem is extremely expensive from a computational point of view. 
In this paper, we exploit the low-rank structure in the solution of the Vlasov equation. More specifically, we consider the Vlasov--Maxwell system and propose a dynamic low-rank integrator. The key idea is to approximate the dynamics of the system by constraining it to a low-rank manifold. This is accomplished by a projection onto the tangent space. There, the dynamics is represented by the low-rank factors, which are determined by solving lower-dimensional partial differential equations. The proposed scheme performs well in numerical experiments and succeeds in capturing the main features of the plasma dynamics. We demonstrate this good behavior for a range of test problems. The coupling of the Vlasov equation with the Maxwell system, however, introduces additional challenges. In particular, the divergence of the electric field resulting from Maxwell's equations is not consistent with the charge density computed from the Vlasov equation. We propose a correction based on Lagrange multipliers which enforces Gauss' law up to machine precision.

\section{Introduction}

A kinetic description of plasmas is required in many applications; for example, in magnetic confined fusion.
In this paper, we are interested in the Vlasov equation 
\[
\partial_tf(t,x,v)+v \cdot \nabla_x f(t,x,v)- F(f) \cdot \nabla_v f(t,x,v) = 0. \\
\]
The Vlasov equation describes a collisionless plasma, where $f$ is the density of particles in $\mathbb{R}^{d_x+d_v}$, with $d_x,d_v=1,2,3$. The physically interesting case is $d_v \geq d_x$.
The particles are subjected to a force $F$, which gives rise to a nonlinear term. 
When considering a plasma interacting with an electromagnetic field the Lorentz force has to be considered. It is given by $F=E+v\times B$, where $E$ and $B$ denote the electric and magnetic fields, respectively. Moreover, the Vlasov equation has to be coupled with Maxwell's equations, which describe the time evolution of the electric and magnetic fields subject to two divergence constraints; Gauss' laws for electricity and magnetism, respectively. 

The numerical solution of the Vlasov--Maxwell equations faces several challenges. Those difficulties are primarily caused by the high dimensionality of the system. A full discretization of phase space with $n_x$ grid points per direction in $x$ and $n_v$ points per direction in $v$ leads to a storage cost of $\mathcal{O}(n_x^{d_x}n_v^{d_v})$, which is extremely expensive. This is often referred to as the curse of dimensionality.
To alleviate this problem, particle in cell (PIC) methods have been widely used. In this approach, only the physical space is discretized, and a conglomerate of particles is advanced by following the characteristic curves of the Vlasov equation; see, e.g., the review paper \cite{verboncoeur2005particle}. 
Numerical integrators based on fully discretizing phase space, i.e. the Eulerian approach, have become affordable only in recent years, both due to improved algorithms and the greatly increased performance of computer systems (see, e.g., \cite{crouseilles2009parallel,rozar2013,einkemmer2015,einkemmer2016mixed,hittinger2013block,wettervik2017relativistic,sonnendrucker1999,vinas2018flux}). 
Algorithms based on splitting methods have been considered. The idea in this setting is to split the problem into a sequence of simple advection equations and then to solve those by semi-Lagrangian or spectral methods. We refer, e.g., to \cite{cheng1976,klimas1994,CEF15,CEF16,sircombe2009,mangeney2002}.

Solvers for the Vlasov--Poisson equations based on low-rank approximations were proposed in the last few years in order to reduce computation and storage cost \cite{Ehrlacher2017,Kormann15}. In both papers, the authors discretize first, perform a time step, and then use tensor truncation to project back to the low-rank approximation space.  
Recently, a different approach was proposed in the context of the Vlasov--Poisson equations; see \cite{EL18}  and the follow-up papers \cite{E18,EL18_cons}. The key idea is to approximate the dynamics of the system by constraining it to the low-rank manifold. The solution is then represented in low-rank form already at the continuous level, before any space or time discretization is performed.  Differential equations for the low-rank factors of the solutions are derived, such that the dynamics of the systems is shifted to the low-rank factors representing the solution. These partial differential equations can then be discretized by any suitable numerical method. A convergence proof for a simplified model, the radiative transfer equation, has been conducted in \cite{ding2009}.

An integrator based on this approach was first considered in \cite{LO14} in the context of matrix differential equations. It is known as the projector-splitting integrator for the dynamic low-rank approximation. More theoretical insights were given in \cite{KL07,KLW16} and extensions for various algorithms and tensor formats can be found in \cite{Koch2010,lubich15tio,lubich13dab,Lubich2017}. 
Moreover, the dynamic low-rank approximation was employed for different types of problems, originally in quantum mechanics; see \cite{meyer90tmc,meyer09mqd,Lubich2008,lubich15tii}.

In the context of kinetic equations, the low-rank integrator has the drawback that it does not preserve the physical structure of the system. That is, physical invariants are in general not preserved leading to severe limitations on the long time integration; see \cite{EL18}. In \cite{EL18_cons} the authors studied this issue and proposed a quasi-conservative scheme for the preservation of mass and momentum for the Vlasov--Poisson equations.

In this paper, we propose a dynamic low-rank integrator for the Vlasov--Maxwell equations. In contrast to the Vlasov--Poisson equations, we need to take into account the influence of the magnetic field.
The obtained differential equations, which describe the dynamics on the low-rank manifold, are advection equations in a lower-dimensional subspace and therefore well suited for applying semi-Lagrangian schemes or Fourier techniques. The proposed method never requires the computation of quantities that lie outside of the low-rank approximation space and therefore has a computational advantage compared to the method in \cite{Kormann15}.

The coupling of the Vlasov equation with Maxwell's equations poses additional challenges. 
In the Vlasov--Poisson equations the values of the density function and the electric field are consistent for the entirety of the time integration process. The value of the charge density obtained from the Vlasov equation is directly used to compute the electric field. In the Vlasov--Maxwell equations, this approach is not possible, since the electric and magnetic fields are themselves coupled. Thus, the charge density obtained from the Vlasov equation is, in general, not consistent with the divergence of the electric field obtained from Maxwell's equations. This leads to a violation of Gauss' law.  
We propose to include a correction strategy in the low-rank algorithm, such that the divergence constraint for the electric field is enforced at each integration step. We follow the ideas in \cite{ADHRS93} and associate a Lagrange multiplier with the constraint on the electric field. 

The paper is structured as follows. In section \ref{sect:V-M} we introduce the Vlasov--Maxwell equations. We then motivate the use of a low-rank approach in section \ref{sec:low-rank-linear}. The proposed dynamic low-rank integrator is explained in detail in section \ref{sect:proj-split}. In section \ref{sect:correction} we explain a correction technique that enforces the divergence constraint at each time step. The description of the setting for the numerical experiments and the implementation details can be found in section \ref{sect:implementation}. In sections \ref{sect:landau}, \ref{sect:two-stream}, \ref{sect:bump}, and \ref{sect:weibel} we illustrate the method by means of some numerical experiments. We consider a Landau-type problem, a two-stream instability, a bump-on-tail instability, and a Weibel instability, respectively. While the basic low-rank integrator succeeds in capturing the main dynamics of the system, we observe that a number of invariants are not preserved and, in particular, Gauss' law is strongly violated. 
The numerical comparison between the basic low-rank integrator and the divergence preserving low-rank integrator is also carried out here.

\section{The Vlasov--Maxwell equations} \label{sect:V-M}
In this section, we introduce the Vlasov--Maxwell equations in detail. We restrict our discussion to one particle species, in our case electrons. However, the extension to multiple species is straightforward.

We denote the phase space by $\Omega_{x}\times\Omega_{v}$ and write $(x,v) \in \Omega_{x}\times\Omega_{v}$. 
We consider here the dimensionless form of the Vlasov equation coupled with Maxwell's equations
\begin{align}
& \partial_tf(t,x,v)+v \cdot \nabla_x f(t,x,v)-\left(E(t,x)+v\times B(t,x)\right) \cdot \nabla_v f(t,x,v) = 0 \label{eq:vlasov}\\
& \partial_t E(t,x) = c^2\, \nabla \times B(t,x)- j(t,x)  \label{eq:ev_E}\\
& \nabla \cdot E(t,x) = n_0- \rho(t,x) \label{eq:Gauss} \\
& \partial_t B(t,x) = -\nabla \times E(t,x) \label{eq:ev_B}\\
& \nabla \cdot B(t,x) = 0.  \label{eq:null_div}
\end{align}
The unknown $f(t,x,v)$ denotes the electron density at time $t$ at the location $(x,v)$, $E$ and $B$ are the electric and magnetic fields, respectively. In the following we set the speed of light $c$ to the value 1.
The charge and current densities are defined as
\[
\rho(t,x) = \int_{\Omega_v} f(t,x,v) \dd v \quad \text{and} \quad j(t,x) = \int_{\Omega_v} v f(t,x,v) \dd v,
\]
respectively.
Since we consider only one particle species here, the value $n_0 = \frac{1}{\lvert \Omega_x \rvert}\int_{\Omega_{x}} \rho(t,x) \dd x$ normalizes the electron charge density with respect to the total charge of the surrounding ions. This corresponds to assuming a uniform background distribution of ions. 

The evolution of the electromagnetic field is restricted by two divergence constraints, the so called Gauss laws for electricity \eqref{eq:Gauss} and magnetism \eqref{eq:null_div}. These equations are automatically satisfied in the continuous case provided that the initial values $E(t=0,x)$ and $B(t=0,x)$ satisfy the same properties. 
In fact, Gauss' law \eqref{eq:Gauss} follows from the continuity equation 
\[ 
\partial_t \rho+ \nabla \cdot j = 0.
\] 
In other words, Gauss' law is a direct consequence of the conservation of charge. 
Further, recalling that $\nabla \cdot (\nabla \times \nobreak \cdot) = 0$, the zero divergence condition \eqref{eq:null_div} for $B$ holds. 
Hence, most Maxwell solvers take into account only the time evolution of the two fields, i.e., \eqref{eq:ev_E} and \eqref{eq:ev_B}, and neglect the divergence constraints. However, one usually has to be careful that these physical constraints are still satisfied in the numerical method. Otherwise, the solution can become inaccurate, and in some situations, even numerical instabilities are known to develop. We refer the reader to \cite{MHD02} and references contained therein.
The zero divergence constraint for the magnetic field can be fulfilled by selecting an appropriate space discretization. The discrete differential operator $\widetilde{\nabla}$, a discretization of $\nabla$, has to be chosen such that $\widetilde{\nabla} \cdot (\widetilde{\nabla} \times \cdot) = 0$. One possibility is to use spectral approximations. A staggered grid can also be employed, where the components of the electric and magnetic fields are considered at different grid points, the so called Yee discretization; see \cite{Yee66}. However, this is not as easily possible for Gauss' law for electricity as the charge density couples to the Vlasov equation.

Solving the Vlasov equation is very expensive due to its high dimensionality, amongst other factors. 
Recently, a low-rank approach was proposed for the Vlasov--Poisson equations; see \cite{EL18,EL18_cons}.
The key assumption is that the solution has low-rank structure, i.e., it can be approximated as a combination of a relatively small number of basis functions that only depend either on $x$ or on $v$. Then, a function $f$ of rank $r$ can be written as
\[
f(t,x,v) = \sum_{i,j=1}^r X_i(t,x) S_{ij}(t) V_j(t,v).
\]
Such an approximation will be employed here.

\section{Low-rank representation of the Vlasov--Maxwell equations \label{sec:low-rank-linear}}
Our goal in this section is to analytically investigate a situation
in which the solution of the Vlasov\textendash Maxwell equations admits
a low-rank representation. This will be helpful for the interpretation
of the numerical results that are obtained in sections  \ref{sect:landau}, \ref{sect:two-stream}, \ref{sect:bump}, and \ref{sect:weibel}. 

We start by linearizing the Vlasov\textendash Maxwell equations around an equilibrium, as done, for example, in \cite{Despres}.
To that end we introduce quantities $f(t,x,v)=f^{(0)}(v)+f^{(1)}(t,x,v)$,
where $f^{(0)}$ is the unperturbed equilibrium, $E(t,x)=0+E^{(1)}(t,x)$,
and $B(t,x)=0+B^{(1)}(t,x)$. This assumes that our initial value
is a perturbation of the equilibrium density given by $f^{(0)}(v)$,
say a Gaussian, and that no electric or magnetic fields are required
to maintain the equilibrium. The linearized Vlasov\textendash Maxwell
equations are then given by 
\begin{align*}
 & \partial_{t}f^{(1)}(t,x,v)+v\cdot\nabla_{x}f^{(1)}(t,x,v)-(E^{(1)}(t,x)+v\times B^{(1)}(t,x))\cdot\nabla_{v}f^{(0)}(v)=0\\
    & \partial_{t}E^{(1)}(t,x)=\nabla\times B^{(1)}(t,x)-j^{(0)}-\int_{\Omega_{v}}vf^{(1)}(t,x,v)\,\mathrm{d}v\\
 & \partial_{t}B^{(1)}(t,x)=-\nabla\times E^{(1)}(t,x),
\end{align*}
where $j^{(0)} = \int_{\Omega_v} v f^{(0)}(v)\,\mathrm{d}v$ is the (constant) current of the equilibrium. Here, we have neglected higher order terms. Then we perform a Fourier transform in $x$ and denote the corresponding
frequencies by $k$. This gives
\begin{align}
 & \partial_{t}\hat{f}_{k}^{(1)}(t,v)+i(v\cdot k)\hat{f}_{k}^{(1)}(t,v)-\left(\hat{E}_{k}^{(1)}(t)+v\times\hat{B}_{k}^{(1)}(t)\right)\cdot\nabla_{v}\hat{f}^{(0)}(v)=0\label{eq:linearized-vlasov-fourier}\\
    & \partial_{t}\hat{E}_{k}^{(1)}(t)=ik\times\hat{B}_{k}^{(1)}(t)- j^{(0)} \delta_{k,0} -\int_{\Omega_{v}}v\hat{f}_{k}^{(1)}(t,v)\,\mathrm{d}v\label{eq:linearized-E-fourier}\\
 & \partial_{t}\hat{B}_{k}^{(1)}(t)=-ik\times\hat{E}_{k}^{(1)}(t).\label{eq:linearized-B-fourier}
\end{align}
Thus, all of the Fourier modes decouple.

Now, let us assume that 
\[
f^{(1)}(0,x,v)=\sum_{i=1}^{m}\hat{f}_{k_{i}}^{(1)}(0,v)e^{ik_{i}x},\qquad B^{(1)}(0,x)=\sum_{i=1}^{m}\hat{B}_{k_{i}}^{(1)}(0)e^{ik_{i}x}.
\]
That is, the initial perturbation of the equilibrium $f^{(0)}(v)$ and of the magnetic field 
can be written as a linear combination of $m$ Fourier modes. This
is not a very restrictive assumption. In fact, most of the initial
values for the well known test problems can be written in this form
(often even for $m=1$).

Now, if $k\neq 0$ and $k\neq k_{i}$ for all $i\in\{1,\dots,m\}$ we have $\hat{f}_{k}^{(1)}(0,v)=0$.
We initialize the electric field using Gauss' law in Fourier space, i.e.
\[
\hat{E}_{k}^{(1)}(0)=\frac{ik}{k^{2}}\int_{\Omega_{v}}\hat{f}_{k}^{(1)}(0,v)\,\mathrm{d}v=0.
\]
Then, the unique solution
of (\ref{eq:linearized-vlasov-fourier})--(\ref{eq:linearized-B-fourier})
is given by
\[
\hat{f}_{k}^{(1)}(t,v)=0,\qquad\hat{E}_{k}^{(1)}(t)=0,\qquad \hat{B}_{k}^{(1)}(t)=0.
\]
Thus, no new Fourier modes are excited (except possibly the constant mode $k=0$) and it follows that
the solution of the linearized Vlasov\textendash Maxwell equations can
be written as
\begin{equation}
f(t,x,v)=f^{(0)}(v)+\sum_{i=0}^{m}\hat{f}_{k_{i}}^{(1)}(t,v)e^{ik_{i}x}\label{eq:solution-linearized}
\end{equation}
with $k_0=0$. This is a low-rank approximation of rank at most $m+1$, which is
the desired result.

The full nonlinear dynamics of the Vlasov\textendash Maxwell equations
can, of course, excite new Fourier modes. If the new modes reach a
significant amplitude, the linear analysis conducted above is no longer
accurate. However, the analysis is still useful since it explains
a behavior that we observe in a number of numerical experiments (see
also sections \ref{sect:landau}, \ref{sect:two-stream}, \ref{sect:bump}, and \ref{sect:weibel}). Namely, that a very small rank is sufficient
to almost exactly represent the dynamics of the Vlasov equation in
the linear regime. This is the case, for example, for the two-stream
instability before the electric energy saturates or for linear Landau
damping. 

It is also interesting to note that the analysis conducted does not
exploit the full potential of the low-rank approximation considered
in this paper. More specifically, the low-rank structure of (\ref{eq:solution-linearized})
assumes that $X_{i}(x)=e^{ik_{i}x}$. In contrast, the low-rank approximation
considered in this paper makes no such stipulation. Consequently, in our algorithm
the $X_{i}$ can be arbitrary functions of $x$. More details will be given in the next section. Hence, the
class of functions that can be represented by a low-rank approximation
is much wider than what the Fourier analysis above would indicate.
Thus, being in the \textit{linear regime} is a sufficient but not necessary condition to obtain a low-rank approximation with small
rank. As an example, we note that the the saturation of the electric
energy (in the two-stream or Weibel instability) is resolved very
well even for a small rank (see sections \ref{sect:two-stream} and \ref{sect:weibel}). Saturation is a nonlinear
phenomenon which is not at all captured by the linear theory (i.e.~by
the dispersion relation).

It is also interesting to note at this point that the solution can
have a very high degree of filamentation, but still be low-rank. This
fact was explained and then used in \cite{EL18}
to efficiently solve the plasma echo problem in the electrostatic
regime. In the above analysis, filamentation is implicitly contained
in functions such as $f_{k_{i}}^{(1)}(t,v)$, which can be oscillatory
(and thus have large derivatives) in $v$. This is a significant advantage of the low-rank
approach, compared to other dimension reduction techniques.

In the next section, we summarize the framework for the construction of a low-rank integrator, we define the low-rank manifold in which the low-rank solution lives and describe its tangent space. We then use these two structures to derive a low-rank integrator.

\section{The low-rank projector-splitting integrator} \label{sect:proj-split}
In this section, we explain the low-rank integrator for the Vlasov--Maxwell equations for the 3+3 dimensional case in full detail, following the ideas in the seminal paper \cite{LO14} and the recent works \cite{EL18,E18,EL18_cons}.

First, we summarize the framework for the construction of the integrator. 
We assume that the particle density function $f$ can be represented in low-rank format in the following way 
\begin{equation}\label{eq:l-r}
f(t,x,v) = \sum_{i,j=1}^r X_i(t,x) S_{ij}(t) V_j(t,v), 
\end{equation}
where $r$ indicates the rank. The functions $X_i(x), \ i = 1, \dots,r$ express the dependence of $f$ on the space variable $x \in \Omega_x$ and the functions $V_j(v), \ j = 1, \dots,r$ the dependence of $f$ on the velocity variable $v \in \Omega_v$. In the following we drop the range of the summation indices and always assume that they run from 1 to $r$.

We seek an approximation of $f$ in the low-rank manifold, defined as follows 
\begin{align*}
\M = \{& f \in L^2(\Omega): f(x,v) = \sum_{i,j} X_i(x) S_{ij} V_j(v) \ \text{with} \ S = (S_{ij}) \in \R^{r \times r}   \text{ invertible}, \\
& \ X_i(x) \in L^2(\Omega_x), \ V_j(v) \in L^2(\Omega_v), \ \langle X_i, X_k\rangle _x = \delta_{ik}, \ \langle V_j, V_l\rangle _v = \delta_{jl} \},
\end{align*}
where $\langle \cdot,  \cdot\rangle _x$ and $\langle \cdot, \cdot \rangle _v$ denote the $L^2$ inner products in $\Omega_x$ and $\Omega_v$, respectively. 
The corresponding tangent space is 
\begin{align*}
\TM = \{ \dot{f} &\in L^2(\Omega): \dot{f}(x,v)=\sum_{i,j} \left(\dot{X}_i(x) S_{ij} V_j(v) + X_i(x) \dot{S}_{ij} V_j(v)+ X_i(x) S_{ij} \dot{V}_j(v)\right) \\ 
& \text{with} \ \dot{S} = (\dot{S}_{ij}) \in \R^{r \times r}, \ \dot{X}_i(x) \in L^2(\Omega_x), \ \dot{V}_j(v) \in L^2(\Omega_v), \ \langle X_i, \dot{X}_k\rangle _x = 0, \ \langle V_i,\dot{V}_j\rangle _v = 0 \}. \end{align*}
We aim to approximate the dynamics of the Vlasov equation \eqref{eq:vlasov} by constraining it to the manifold $\M$. This is done by means of an orthogonal projection onto the tangent space $\TM$. 
The projector $P(f)$ is given by the following expression 
\[ P(f) g = \sum_j \langle V_j,g\rangle _v V_j - \sum_{i,j} X_i \langle X_iV_j, g\rangle _{x,v} V_j + \sum_i X_i \langle X_i,g\rangle _x. 
\]
By introducing the following two vector spaces $\bar{X} = \text{span}\{X_1, \dots, X_r\}$ and $\bar{V} = \text{span}\{V_1, \dots, V_r\}$,
the projection can be written as 
\begin{equation} \label{eq:projection} 
P(f)g = P_{\bar{V}} g - P_{\bar{V}} P_{\bar{X}} g + P_{\bar{X}} g.
\end{equation}

In \cite{LO14} the authors proposed a splitting integrator for solving the projected equation
\[ 
\partial_t f = - P(f) \left(v \cdot \nabla_x f-(E+v\times B) \cdot \nabla_v f \right). 
\]
We split the right hand side into three terms and solve the arising subproblems:
\begin{align} 
&\partial_t f = -P_{\bar{V}} \left(v \cdot \nabla_x f-(E+v\times B) \cdot \nabla_v f \right)  \label{eq:Keq} \\
 &\partial_t f = P_{\bar{V}} P_{\bar{X}}\left(v \cdot \nabla_x f-(E+v\times B) \cdot \nabla_v f \right)  \label{eq:Seq} \\ 
 &\partial_t f = -P_{\bar{X}} \left(v \cdot \nabla_x f-(E+v\times B) \cdot \nabla_v f \right). \label{eq:Leq} 
\end{align} 

In the following we explain in detail the first-order scheme, the so called Lie splitting. The second-order integrator is described in algorithm \ref{alg:Strang}.
At first we consider the projector $P_{\bar{V}}$ and look at subproblem \eqref{eq:Keq}. 
We represent the approximate particle density $f$ as 
\[ 
f(t,x,v) = \sum_{j} K_j(t,x) V_j(t,v) \quad \text{with} \quad K_j(t,x)= \sum_{i} X_i(t,x) S_{ij}(t).
\]
The initial value is represented by $X_i^0(x)=X_i(0,x)$, $S_{ij}^0=S_{ij}(0)$ and $V_j^0(v)=V_j(0,v)$.
Then, equation \eqref{eq:Keq} can be rewritten as 
\begin{align*} 
& \sum_j \partial_t K_j(t,x)V_j(t,v) + \sum_j K_j(t,x) \partial_t V_j(t,v) = \\
& -\sum_j \langle V_j(t,\cdot),\left[v\mapsto v \cdot \nabla_x f(t,x,v)- \left( E(x)+ v \times B(x) \right) \cdot \nabla_v f(t,x,v)\right]\rangle _v V_j(t,v). 
\end{align*}
By making use of the orthogonality conditions $\langle V_j, V_l\rangle _v = \delta_{jl}$ and $\langle V_i,\dot{V}_j\rangle _v = 0$ we obtain the following differential equation for the coefficients in the basis expansion 
\begin{equation} \label{eq:Kstep}
\partial_t K_j(t,x) = -\sum_l c_{jl}^1 \cdot \nabla_x K_l(t,x) + \sum_l c_{jl}^2 \cdot E(t,x) K_l(t,x) + \sum_l c_{jl}^3 \cdot B(t,x) K_l(t,x).
\end{equation}
Moreover, from the previous calculation it also follows that $V_j(t,v) = V_j^0(v)$.  
The coefficients $c_{jl}^{1,2} \in \R^3$ are given by the following expressions
\[c_{jl}^1 = \int_{\Omega_v} v V_j^0(v) V_l^0(v) \dd v, \quad c_{jl}^2 = \int_{\Omega_v} V_j^0(v) \nabla_v V_l^0(v) \dd v.
\]
The last term in \eqref{eq:Kstep} was rewritten by means of the relation between scalar and vector products $(a \times b) \cdot c = (c\times a) \cdot b$. The coefficient $c_{jl}^3 \in \R^3$ is 
\[ c_{jl}^3 = \int_{\Omega_v} V_j^0(v) \left(\nabla_v V_l^0(v) \times v\right) \dd v.\]
Thus, we integrate equation \eqref{eq:Kstep} with initial value $K_j(0,x) = \sum_{i} X_i^0(x) S_{ij}^0$ to obtain $K_j(\tau,x) = K_j^1(x)$, where $\tau$ denotes the integration step. After performing a QR-decomposition
\[
K_j^1(x) = \sum_{i} X_i^1(x) S_{ij}^1,
\] 
we obtain orthonormal $X_j^1(x)$, and $S_{ij}^1$.

The subproblem \eqref{eq:Seq} is treated in a similar way. In this step the particle function is represented as $f(t,x,v) = \sum_{i,j} X_i(t,x)S_{ij}(t)V_j(t,v)$ and only $S_{ij}$ is updated. The factors $X_i(t,x)$ and $V_j(t,v)$ are unchanged, i.e., $X_i(t,x)=X_i^1(x)$ and $V_j(t,v)= V_j^0(v)$. The corresponding differential equation is 
\begin{equation} \label{eq:Sstep} 
\dot{S}_{ij}(t) = \sum_{k,l} \left(c_{jl}^1 \cdot d_{ik}^2 - c_{jl}^2 \cdot d_{ik}^1[E(t,x)]-c_{jl}^3  \cdot d_{kl}^3[B(t,x)]\right) S_{kl}(t),
\end{equation}
which has to be integrated with initial value $S_{ij}(0) = S_{ij}^1$. After one integration step we obtain $S_{ij}( \tau) = S_{ij}^2$.
The coefficient $d_{ik}^2 \in \R^3$ is given by 
\[ d_{ik}^2 = \int_{\Omega_{x}} X_i^1(x) \nabla_xX_k^1(x) \dd x, 
\]  
whereas the coefficients $d_{ik}^{1,3} \in \R^3$ are defined as follows
\[ d_{ik}^1[E](t,x) = \int_{\Omega_{x}} X_i^1(x) E(t,x) X_k^1(x) \dd x, \quad d^3_{ik}[B](t,x) = \int_{\Omega_{x}} X_i^1(x) B(t,x)  X_k^1(x) \dd x.
\]
For the solution of \eqref{eq:Leq} we represent the particle function as 
\[
f(t,x,v) = \sum_{i} X_i(t,x) L_i(t,v)  \quad \text{with} \quad L_i(t,x)= \sum_{j} S_{ij}(t) V_j(t,v).
\]
The corresponding differential equation is 
\begin{equation} \label{eq:Lstep}
\partial_t L_i(t,v) = -\sum_k d_{ik}^2 \cdot v L_k(t,v) + \sum_k d_{ik}^1[E(t,x)] \cdot \nabla_v L_k(t,v) + \sum_k d_{ik}^3[B(t,x)] \cdot (\nabla_v L_k(t,v) \times v).
\end{equation}
We integrate it with initial value $L_i(0,x)= \sum_{j} S^2_{ij} V_j^0(v)$ and obtain $L_i(\tau,v) = L_i^1(v)$. After a QR-decomposition we get orthonormal $V_j^1(v)$, and $S_{ij}^3$ such that 
\[ 
L_i^1(v)= \sum_{j} S^3_{ij} V_j^1(v).
\] 
The value of the particle function after one time step $f(\tau,x,v)$ is approximated by
\[  
f^1(x,v) = \sum_{i,j} X_j^1(x) S_{ij}^3 V_j^1(v).
\]
Note that the three differential equations \eqref{eq:Kstep}, \eqref{eq:Sstep} and \eqref{eq:Lstep} are solved by a favorable numerical scheme, where the electric and magnetic fields are evaluated at the beginning of the time step. 
 
Further, the solution of the evolution equations for the electric and magnetic fields, \eqref{eq:ev_E} and \eqref{eq:ev_B} respectively, has to be computed. We remark that the space discretization has to be chosen such that the zero divergence constraint \eqref{eq:null_div} is satisfied. We suggest either a spectral discretization or the Yee space discretization, see \cite{Yee66}, and denote the discretization of the operator $\nabla$ by $\widetilde{\nabla}$. A first-order scheme is sufficient for the time integration, for example
\begin{align}
E(\tau,\cdot) & \approx E^1 =  E^0 + \tau \widetilde{\nabla} \times B^0 - \tau j^0 \label{eq:E_first}\\
B(\tau,\cdot) & \approx B^1= B^0-\tau \widetilde{\nabla} \times E^1 \label{eq:B_first}, 
\end{align}
where $E(0,\cdot)=E^0$, $B(0,\cdot)=B^0$, and the current density is computed as follows 
\[
j(0,x) \approx j^0(x) = \sum_{i,j} X_i^0(x) S_{ij}^0 \left(\int_{\Omega_{v}} v V_j^0(v) \dd v \right).
\]
In algorithm \ref{alg:Strang} we illustrate a second-order scheme. The solution of the Vlasov equation is approximated by means of the Strang splitting scheme. The solution of Maxwell's equations is computed on a staggered grid in time; see e.g.~\cite{Yee66}. 

\begin{algorithm}
	\caption{A second-order low-rank integrator for the Vlasov--Maxwell equations} \label{alg:Strang}
	\begin{algorithmic}[1]
		\Statex \textbf{Input}: $X_i^0$, $S_{ij}^0$, $V_j^0$ such that $f(0) \approx \sum_{i,j} X_i^0 S_{ij}^0 V_j^0$, and $E(0,\cdot)= E^0$, $B( \frac{\tau}{2},\cdot)=B^{1/2}$.
		\Statex \textbf{Output}: $X_i^1$, $S_{ij}^3$, $V_j^1$ such that $f(\tau) \approx \sum_{i,j} X_i^1 S_{ij}^3 V_j^1$, and $E^1 \approx E(\tau,\cdot)$, $B^{3/2} \approx B(\frac{3}{2}\tau,\cdot)$.
		\State Solve equation \eqref{eq:ev_E} up to time $\tau/2$ to get $E^{1/2} = E^0 +\frac{\tau}{2} \widetilde{\nabla} \times B^{1/2} -\frac{\tau}{2}  j^0$. 
		\State Solve equation \eqref{eq:Kstep} with initial value $\sum_{i} X_i^0 S_{ij}^0$ and $E(t,\cdot)=E^{1/2}$, $ B(t,\cdot)=B^{1/2}$ up to time $\tau/2$ to obtain $K^{1/2}_j$.
		\State Orthogonalize $K_j^{1/2}$ by a QR-decomposition to obtain $X_i^{1/2}$ and $S_{ij}^{1/2}$. 
		\State Solve equation \eqref{eq:Sstep} with initial value $S_{ij}^{1/2}$ and $E(t,\cdot)=E^{1/2}$, $ B(t,\cdot)=B^{1/2}$ up to time $\tau/2$ to obtain $S_{ij}^1$.
		\State Solve equation \eqref{eq:Lstep} with initial value $ \sum_{j} S_{ij}^1 V_j^0$ and $E(t,\cdot)=E^{1/2}$, $ B(t,\cdot)=B^{1/2}$ up to time $\tau/2$ to obtain $L_i^{1/2}$.
		\State Compute $j^{1/2}$ using $X_i^{1/2}$ and $L_i^{1/2}$.
		\State Solve equation \eqref{eq:Lstep} with initial value $ \sum_{j} S_{ij}^1 V_j^0$ and $E(t,\cdot)=E^{1/2}$, $B(t,\cdot)=B^{1/2}$ up to time $\tau$ to obtain $L_i^1$.
		\State Orthogonalize $L_i^1$ by a QR-decomposition to obtain $S_{ij}^2$ and $V_j^1$.
		\State Solve equation \eqref{eq:Sstep} with initial value $S_{ij}^2$ and $E(t,\cdot)=E^{1/2}$, $ B(t,\cdot)=B^{1/2}$ up to time $\tau/2$ to obtain $S_{ij}^{5/2}$.
		\State Solve equation \eqref{eq:Kstep} with initial value $\sum_{i} X_i^{1/2} S_{ij}^{5/2}$ and $E(t,\cdot)=E^{1/2}$, $ B(t,\cdot)=B^{1/2}$ up to time $\tau/2$ to obtain $K^1_j$.
		\State Orthogonalize $K_j^1$ by a QR-decomposition to obtain $X_i^1$ and $S_{ij}^3$. 
		\State Solve equation \eqref{eq:ev_E} up to time $\tau$ to get $E^1 =  E^0 + \tau \widetilde{\nabla} \times B^{1/2} - \tau j^{1/2}$. 
		\State Solve equation \eqref{eq:ev_B} up to time $\tau$ to get $B^{3/2} = B^{1/2} - \tau \widetilde{\nabla} \times E^1$. 
	\end{algorithmic}
\end{algorithm}

\section{A divergence correction technique} \label{sect:correction}
The low-rank integrator illustrated in the previous section does not explicitly take into account Gauss' law~ \eqref{eq:Gauss}. Such an approach is motivated by the fact that this divergence constraint is automatically satisfied in the continuous case. This follows immediately from the continuity equation, as already explained in section~\ref{sect:V-M}.
However, at the discrete level the charge is not conserved and the divergence constraint is not satisfied. 
Thus, we propose to explicitly include Gauss' law in the numerical scheme. 

We employ an approach based on Lagrange multipliers, following the ideas in \cite{ADHRS93}. The evolution of the electric field is then described by 
\begin{align}
& \partial_t E = \nabla \times B - j -\nabla \phi \label{eq:ev_E_mod}\\
& \nabla \cdot E = \frac{1}{\lvert \Omega_x \rvert}\int_{\Omega_x} \rho \dd x - \rho \label{eq:Gauss_rho}.
\end{align}
The function $\nabla \phi$  is a correction potential, which is interpreted as the Lagrange multiplier. 
This approach turns out to be a particular case of a more general idea proposed in, e.g., \cite{MOSSV00, MHD02}. There, the evolution equation for the electric field \eqref{eq:ev_E_mod} is coupled to the corresponding constraint via the function $\phi$. The Gauss law is modified to
\[
\mathcal{D} (\phi) + \nabla \cdot E = \frac{1}{\lvert \Omega_x \rvert}\int_{\Omega_x} \rho \dd x -\rho, \]
where $\mathcal{D}(\phi)$ is a differential operator. Depending on the choice of such an operator, different techniques arise. In particular, if $\mathcal{D}(\phi) \equiv 0$ the above system is obtained.

The numerical solution of \eqref{eq:ev_E_mod} coupled with \eqref{eq:Gauss_rho} is computed by an operator splitting technique. We first solve the subflow
\begin{equation} 
\partial_t E = \nabla \times B - j.  \label{eq:ev_E_sub}
\end{equation}
This is equivalent to the original equation and thus any suitable numerical method can be applied.
We explain here only a second-order scheme where we make use of a staggered grid in time; see \cite{Yee66}. As before we denote the discretization of $\nabla$ by $\widetilde{\nabla}$, such that $\widetilde{\nabla}\cdot(\widetilde{\nabla}\times \cdot)=0$.
We compute first $\bar{E}^{n+1}$ as the solution of \eqref{eq:ev_E_sub} after one time step of size $\tau$ and get 
\begin{equation} \label{eq:E_bar}
\bar{E}^{n+1} = E^n + \tau \widetilde{\nabla} \times B^{n+1/2} - \tau j^{n+1/2}. 
\end{equation}
The intermediate value of the electric field $\bar{E}^{n+1}$ is then updated by solving the second subflow
\begin{align}
\partial_t E & = -\nabla \phi
\end{align}
coupled to Gauss' law. We obtain the equation 
\begin{equation} 
E^{n+1} =  \bar{E}^{n+1} - \tau \widetilde{\nabla} \phi  \label{eq:discr_E}
\end{equation}
and
\begin{equation} 
\widetilde{\nabla} \cdot E^{n+1} = \frac{1}{\lvert \Omega_x \rvert}\int_{\Omega_x} \rho^{n+1} \dd x - \rho^{n+1} . \label{eq:discr_div}
\end{equation}
Next, we can determine $\phi$ by solving the Poisson equation
\begin{equation} \label{eq:Poisson}
-\widetilde{\Delta} \phi = \frac{\frac{1}{\lvert \Omega_x \rvert}\int_{\Omega_x} \rho^{n+1} \dd x - \rho^{n+1} - \widetilde{\nabla} \cdot E^n}{\tau} + \widetilde{\nabla} \cdot j^{n+1/2}.
\end{equation}
This equation is obtained by applying the discrete divergence to \eqref{eq:discr_E} and substituting the discrete divergence condition \eqref{eq:discr_div}. Then, the electric field is updated as in \eqref{eq:discr_E}. 
This correction step for the electric field has to be incorporated in the low-rank scheme described in algorithm \ref{alg:Strang}. The resulting second-order integrator is summarized in algorithm \ref{alg:Strang_corr}.

Note that the correction potential $\phi$ is constructed by taking into account the specific numerical integrator we use for the solution of \eqref{eq:ev_E_sub}. The error in Gauss' law made in \eqref{eq:E_bar} is then corrected. 
It is crucial to consider the numerical value of the mass instead of the analytic one in the discrete divergence constraint \eqref{eq:discr_div}, since the employed low-rank integrator is not mass preserving.  
Therefore, the divergence constraint has to be determined at the discrete level in terms of the current numerical values of mass and charge density. 

Note that the divergence correction procedure and the quasi-conservative scheme proposed in \cite{EL18_cons} are complementary approaches. The conservation of mass is enforced on the distribution function, whereas the correction proposed here only acts on the electric field; the distribution function only enters as an input via the charge density and the current density. Thus, both procedures can be combined. 

\begin{algorithm}
	\caption{A second-order low-rank integrator for the Vlasov--Maxwell equations with a divergence correction} \label{alg:Strang_corr}
	\begin{algorithmic}[1]
		\Statex \textbf{Input}: $X_i^0$, $S_{ij}^0$, $V_j^0$ such that $f(0,x,v) \approx \sum_{i,j} X_i^0(x) S_{ij}^0 V_j^0(v)$, and $E(0,\cdot) = E^0 $, $B(\frac{\tau}{2},\cdot) = B^{1/2}$.
		\Statex \textbf{Output}: $X_i^1$, $S_{ij}^3$, $V_j^1$ such that $f(\tau,x,v) \approx \sum_{i,j} X_i^1(x) S_{ij}^3V_j^1(v)$, and $E^1\approx E(\tau,\cdot)$, $B^{3/2} \approx B(\frac{3}{2}\tau,\cdot)$.
		\State Perform steps 1--12 of algorithm \ref{alg:Strang} to obtain $X_i^1, S_{ij}^3$ and $V_i^1$.
		\State Compute $\rho^{1}$ using $X_i^1, S_{ij}^3$ and $V_i^1$.
		\State Solve equation \eqref{eq:ev_E_sub} up to time $\tau$ to get $\bar{E}^1 =  E^0 + \tau \widetilde{\nabla} \times B^{1/2} - \tau j^{1/2}$. 
		\State Solve equation \eqref{eq:Poisson} to get $\phi$ and compute $\widetilde{\nabla}\phi$.
		\State Update the electric field as in \eqref{eq:discr_E} to obtain $E^1$.
		\State Solve equation \eqref{eq:ev_B} up to time $\tau$ to get $B^{3/2} = B^{1/2} - \tau \widetilde{\nabla} \times E^1$. 
	\end{algorithmic}
\end{algorithm}

\section{Implementation details} \label{sect:implementation}

In the numerical experiments reported in the following sections, we consider the reduced model obtained by taking into account one dimension in the physical space and two in the velocity space. We describe this model here. To do that, we write $v = (v_1,v_2)$ and denote by $v_1$ and $v_2$ the two components of the velocity. The phase space is then $\Omega_{x_1} \times \Omega_{v_1} \times \Omega_{v_2}$. The electric field $E = E(t,x_1) = (E_1,E_2)(t,x_1)$ is two dimensional, where $E_1(t,x_1)$ and $E_2(t,x_1)$ denote its components. The magnetic field is orthogonal to the electric field and is one dimensional. We thus have  
\[ E(t,x_1) = \begin{pmatrix}
E_1(t,x_1)  \\ E_2(t,x_1) \\ 0 
\end{pmatrix} 
\qquad \text{and} \qquad 
B(t,x_1) = \begin{pmatrix}
0 \\ 0 \\ B_3(t,x_1) 
\end{pmatrix} .
\]
The Vlasov equation reduces to 
\begin{equation} 
\partial_tf+v_1 \partial_{x_1} f-E \cdot \nabla_v f-B_3 \, \mathcal{J}v \cdot \nabla_vf = 0, \label{eq:1+2Vlasov}
\end{equation}
where the matrix $\mathcal{J}$ is given by
\[\mathcal{J} = \begin{pmatrix}
0 & 1 \\ -1 & 0
\end{pmatrix}.
\]
Further, the Maxwell system is given by
\begin{align} 
& \partial_t E_1(t,x_1) = -\int_{\Omega_{v}} v_1 f(t,x_1,v) \dd v \\
& \partial_t E_2(t,x_1) = -\partial_{x_1} B_3(t,x_1) - \int_{\Omega_{v}} v_2 f(t,x_1,v) \dd v\\
& \partial_{x_1} E_1(t,x_1) = \frac{1}{\lvert \Omega_{x_1} \rvert}\int_{\Omega_{x_1}} \rho(t,x) \dd x -\rho(t,x_1) \label{eq:1+2Poisson} \\ 
& \partial_t B_3(t,x_1) = -\partial_{x_1} E_2(t,x_1). 
\end{align}
Note that the zero divergence condition for the magnetic field is automatically fulfilled.

The Vlasov--Maxwell systems conserves the total energy. In the next sections, we study the behaviour of the proposed integrator with respect to this matter. In this respect, we recall here that the electric and magnetic energies are the energies stored in the electric and magnetic fields, respectively. They are defined as 
\[ 
\mathcal{E}_{e}(t) = \frac{1}{2} \int_{\Omega_{x_1}} \lvert E(t,x_1) \rvert ^2 \dd x_1 \quad  \text{and} \quad \mathcal{E}_{m}(t) = \frac{1}{2} \int_{\Omega_{x_1}} \lvert B(t,x_1) \rvert ^2 \dd x_1.
\] 
The kinetic energy is given by
\[
\mathcal{E}_{k}(t) = \frac{1}{2} \int_{\Omega_{x_1}} \int_{\Omega_v} \lvert v \rvert ^2 f(t,x_1,v) \dd v \dd x_1.
\]

The projector-splitting integrator illustrated in section \ref{sect:proj-split} is slightly simpler in the present case. The coefficients $c_{jl}^{1,3}$ are real numbers, whereas $c_{jl}^2$ is a vector in $\R^2$. Similarly, $d_{ik}^{1,3}$ are real numbers and $d_{ik}^2$ is in $\R^2$. The equations for the substeps can then be easily obtained from \eqref{eq:Kstep}, \eqref{eq:Sstep}, and \eqref{eq:Lstep}.

The numerical discretization of the spatial and velocity operators is carried out in spectral space.
We use the FFTW library \cite{FFTW} for computing the discrete Fourier transform in one dimension. Given $n$ real values $X_k$, $k=0, \dots, n-1$ the discrete Fourier transform is defined as follows
\[ 
\widehat{X}_k = \sum_{j=0}^{n-1}X_j\ee^{-2\pi  i j k/n}, 
\]
where $\widehat{X}_k = \widehat{X}_{n-k}^*$ with $*$ denoting the complex conjugate. Due to this symmetry, only half of the output data has to be stored in computer memory. Moreover, it follows that $\hat{X}_0$ and $\hat{X}_{n/2}$, for $n$ even, are real numbers. 
The imaginary part of these numbers is not stored. This introduces an additional numerical error in the solution of the Poisson equation as we explain in the following.

For the time integration of \eqref{eq:Kstep}, \eqref{eq:Sstep}, and \eqref{eq:Lstep} we employ the \texttt{boost::odeint} library. The DOPRI5 (Dormand--Prince) method of order $5$ is used. In order to keep the error small and to guarantee stability, we use $5$ substeps for each step of the splitting scheme.

For computing the correction to the electric field the Poisson equation \eqref{eq:Poisson} has to be solved first.
In spectral space the solution is given by
\[ 
\widehat{\phi}_k = -\frac{1}{k^2} \widehat{g}_k, \quad k=1,\dots,\lfloor n/2 \rfloor \quad \text{and} \quad \hat{\phi}_0 = 0,
\]
where $g$ denotes the right-hand side of \eqref{eq:Poisson}. 
Then the space derivative of $\phi$ has to be computed. In spectral space it is done by multiplying $\hat{\phi}_k$ by $ik$, i.e.
\[
    (\widehat{ \partial_{x} \phi})_k = -\frac{i}{k} \hat{\phi}_k, \quad k\neq0.
\]
Therefore, for $n$ even, the entry $\widehat{\phi}_{n/2}$ is a purely imaginary number. However, due to the symmetry this information can not be stored in memory and thus a numerical error is introduced by the inverse transformation. 
Normally these errors can be neglected. Here, however, the errors grow as the solution evolves in time. 
An easy remedy is to use an odd number of grid points in the real space $\Omega_x$. We will do this in all the simulations conducted in the following sections.  

Let us discuss the computational complexity and storage requirements of the proposed integrator. If $r$ denotes the chosen rank, then for each time step we have to solve
\begin{itemize}
	\item a system of $r$  advection equations in $x$ (this requires the computation of at most $3d_vr^2$ integrals over $\Omega_v$),
	\item a system of $r^2$ ordinary differential equations (this requires the computation of at most $3d_vr^2$ integrals over $\Omega_x$),
	\item a system of $r$ advection equations in $v$.
\end{itemize}
Additionally, one has to consider the cost for solving Maxwell's equations. These are equations posed in $x$ only. In case of the divergence correction the cost of solving the Poisson equation has to be counted also. However, since we do this in Fourier space this constitutes only a  negligible amount of the total run time of the algorithm. If $n_x$ and $n_v$ degrees of freedom are used in each coordinate direction in $x$ and $v$, respectively, then the storage cost is reduced from $\mathcal{O}(n_x^{d_x}n_v^{d_v})$ to $\mathcal{O}(r(n_x^{d_x}+n_v^{d_v}))$. 
As we are using spectral methods, the cost of solving the evolution equations is proportional up to a logarithm to the number of degrees of freedom. In this case our low-rank algorithm
requires $\mathcal{O}(r(n_x^{d_x}+n_v^{d_v}) )$ arithmetic operations for the evolution equations and $\mathcal{O}(r^2(n_x^{d_x}+n_v^{d_v}))$ operations for computing the integrals.

\section{Numerical results for the Landau-type problem} \label{sect:landau}
We consider first a Landau-type problem (in the sense that this configuration shows Landau damping when the speed of light goes to infinity). This is the same setup as in \cite{CEF16}.
We consider the following initial value
\[ f(0,x_1,v_1,v_2) = \frac{1}{2\pi} \ee^{-\frac{1}{2}(v_1^2+v_2^2)}\left(1+\alpha \cos(kx_1)\right),
\]
where we choose $k = 0.4$ and $\alpha=0.01$. The electric field is initialized by solving Gauss' law \eqref{eq:1+2Poisson}. This yields 
\[
E_1(0,x_1) = -\frac{\alpha}{k} \sin(kx_1), \qquad E_2(0,x_1) = 0,
\]
and the magnetic field at $t=0$ is chosen as 
\[
B_3(0,x_1) = \frac{\alpha}{k} \sin(kx_1).
\]
The space domain here is $\Omega_{x_1} = [0,2\pi/k]$ and the velocity domain is chosen as $\Omega_{v_1}\times \Omega_{v_2} = [-5,5]^2$. 

\begin{figure}
	\includegraphics[width = .5\textwidth]{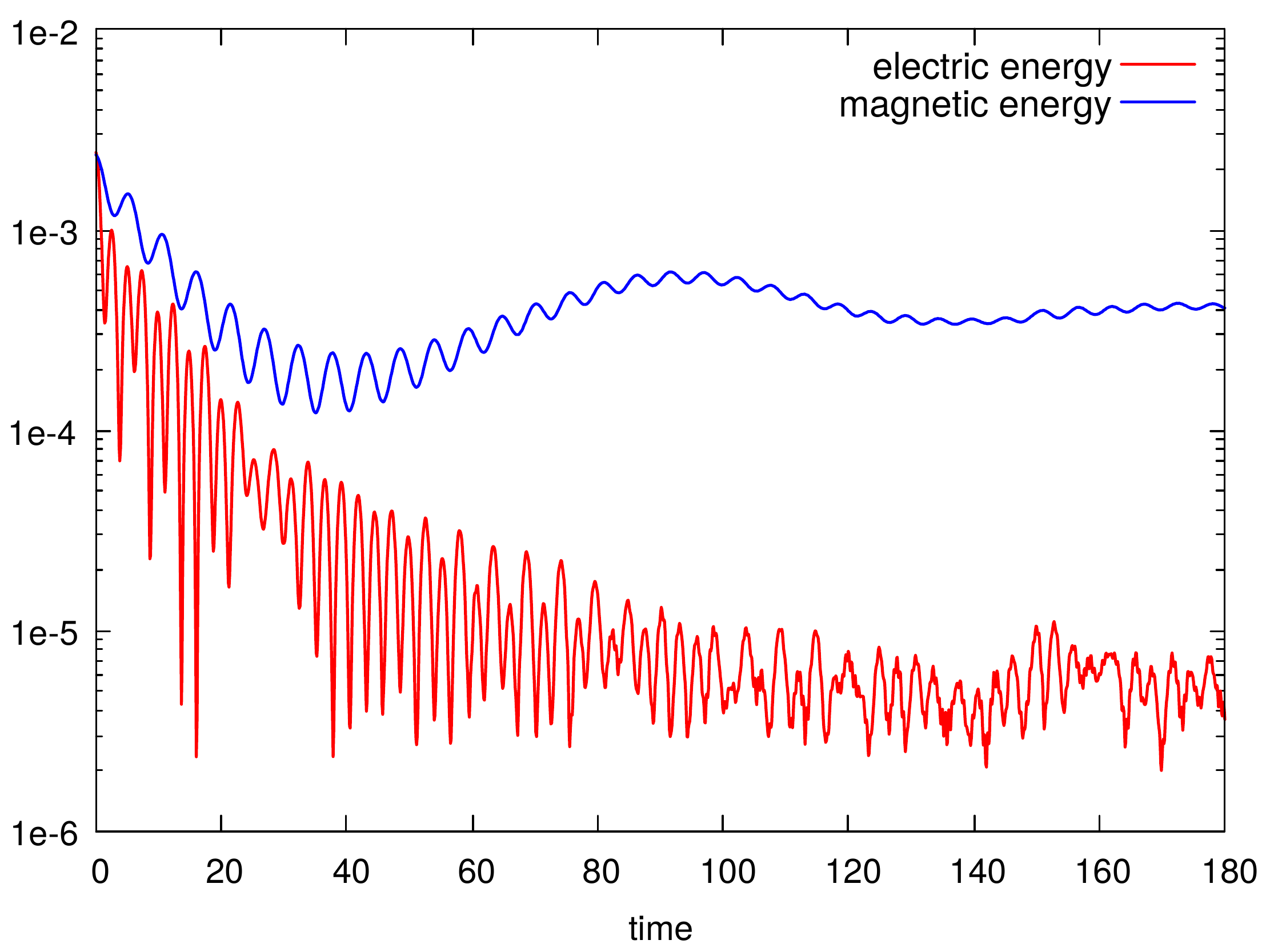}
	\includegraphics[width = .5\textwidth]{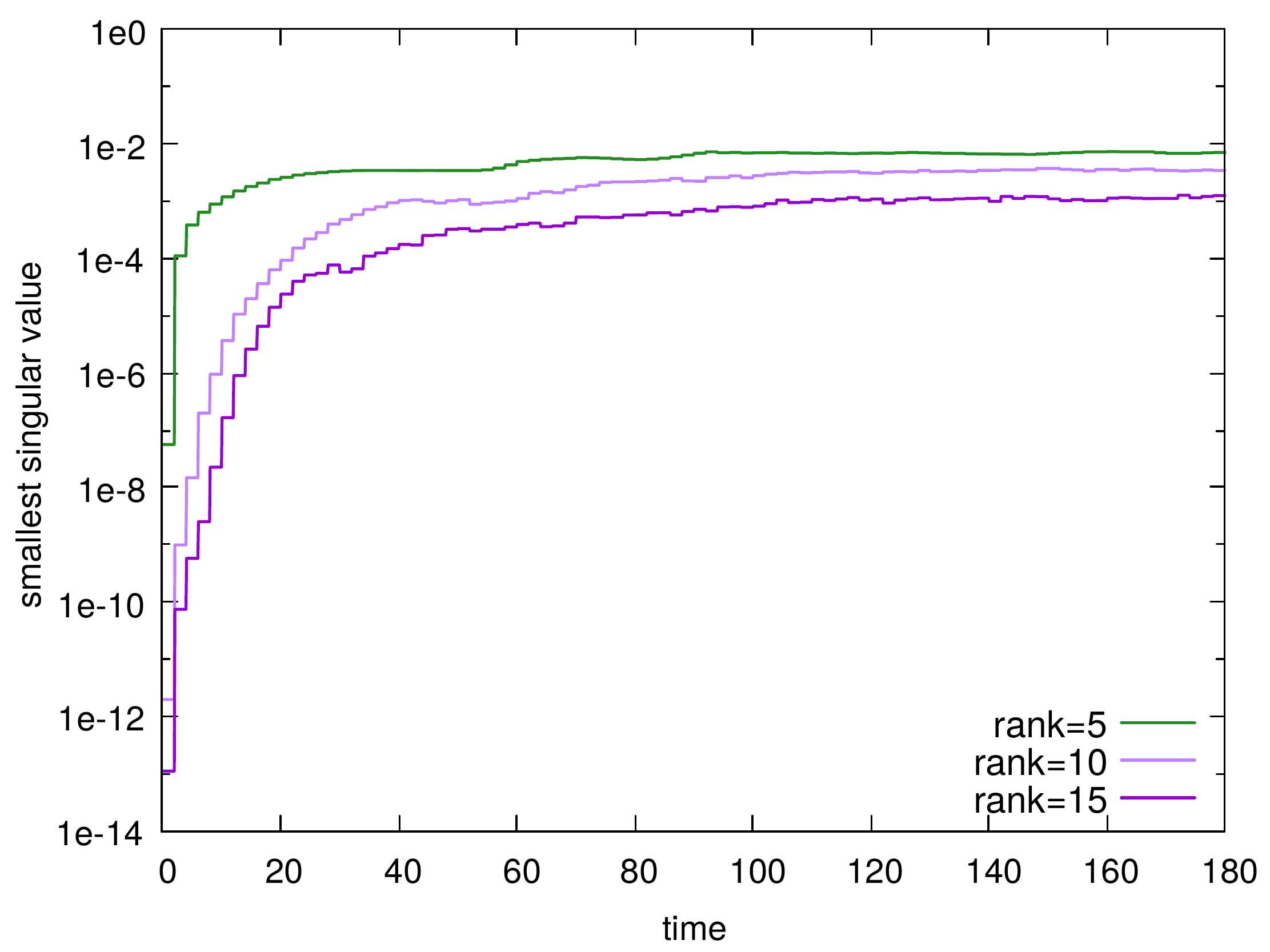} \\
	\center{
	\includegraphics[width = .5\textwidth]{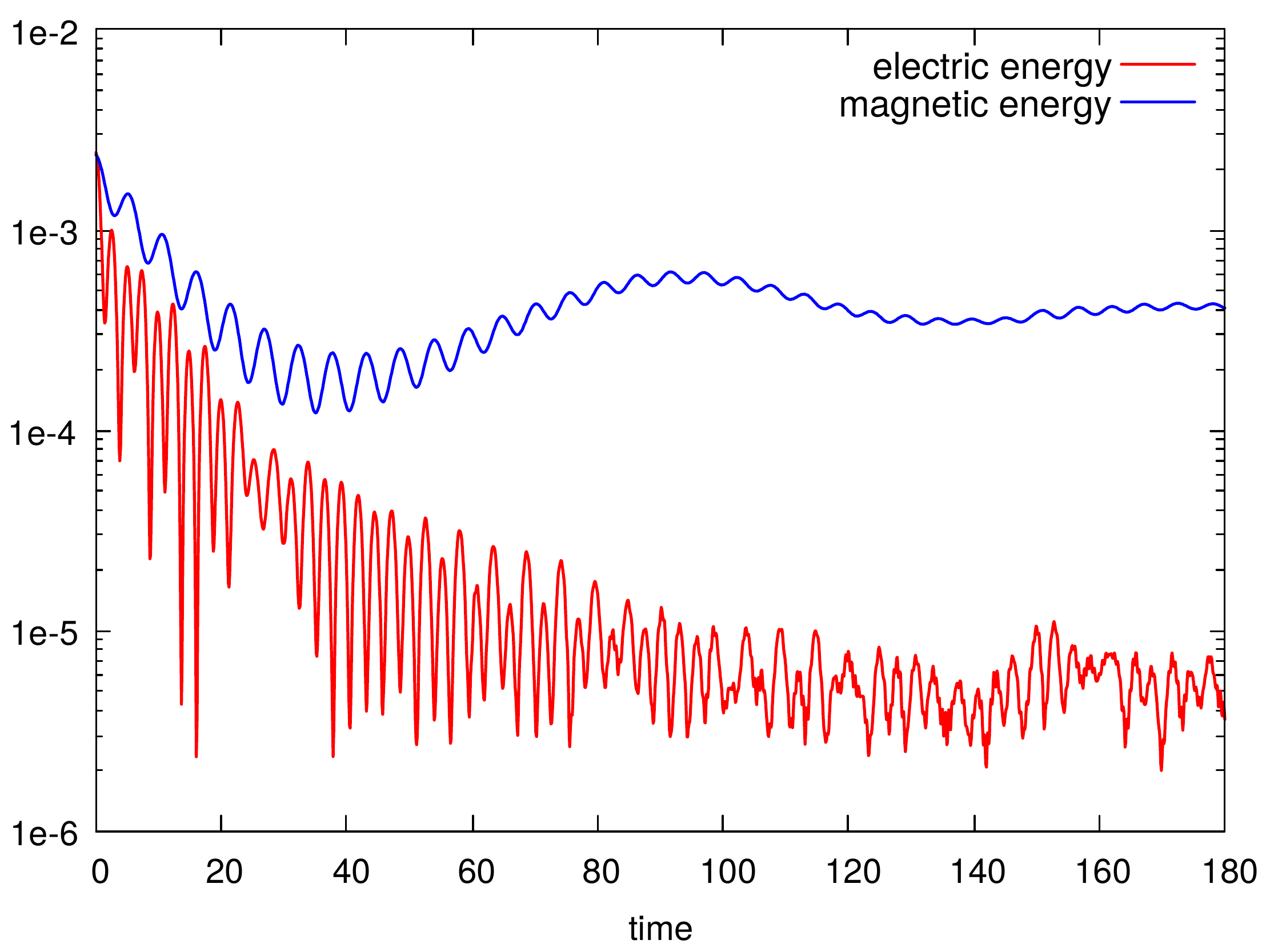}} 
	\caption{Landau-type problem. Results for $n_{x_1}=33$, $n_{v_1}=128$, $n_{v_2}=128$, and $\tau=0.1$ (semi-log scale). The figure on the top left shows the energies obtained with algorithm \ref{alg:Strang} (without correction) for rank = 15; the one on the right displays the smallest singular value of the corresponding distribution function. The bottom figure shows the energies obtained with algorithm \ref{alg:Strang_corr} (with correction) for rank = 15.}
	\label{fig:lt_wc_energy}
\end{figure}

In figure \ref{fig:lt_wc_energy}, on the top, we show the results obtained with algorithm \ref{alg:Strang} (without correction). On the left, we plot the values of the electric and magnetic energies for rank = 15. Our results are very similar to those obtained with an Eulerian Vlasov solver in \cite{CEF16}.
On the right we plot the time evolution of the smallest singular value of the distribution function $f$. The smallest singular value of the matrix $S = (S_{ij})$ gives a measure of the error introduced by truncating the rank. We observe that this error stays bounded in time for all the approximation ranks we tested. Moreover, the choice of a higher approximation rank gives better results, as expected. 

In the third plot in figure \ref{fig:lt_wc_energy} we display the time evolution of the electric and magnetic energies for the divergence preserving scheme given in algorithm \ref{alg:Strang_corr}. We observe that the results are very close to the uncorrected results.

\begin{figure}
	\includegraphics[width = .5\textwidth]{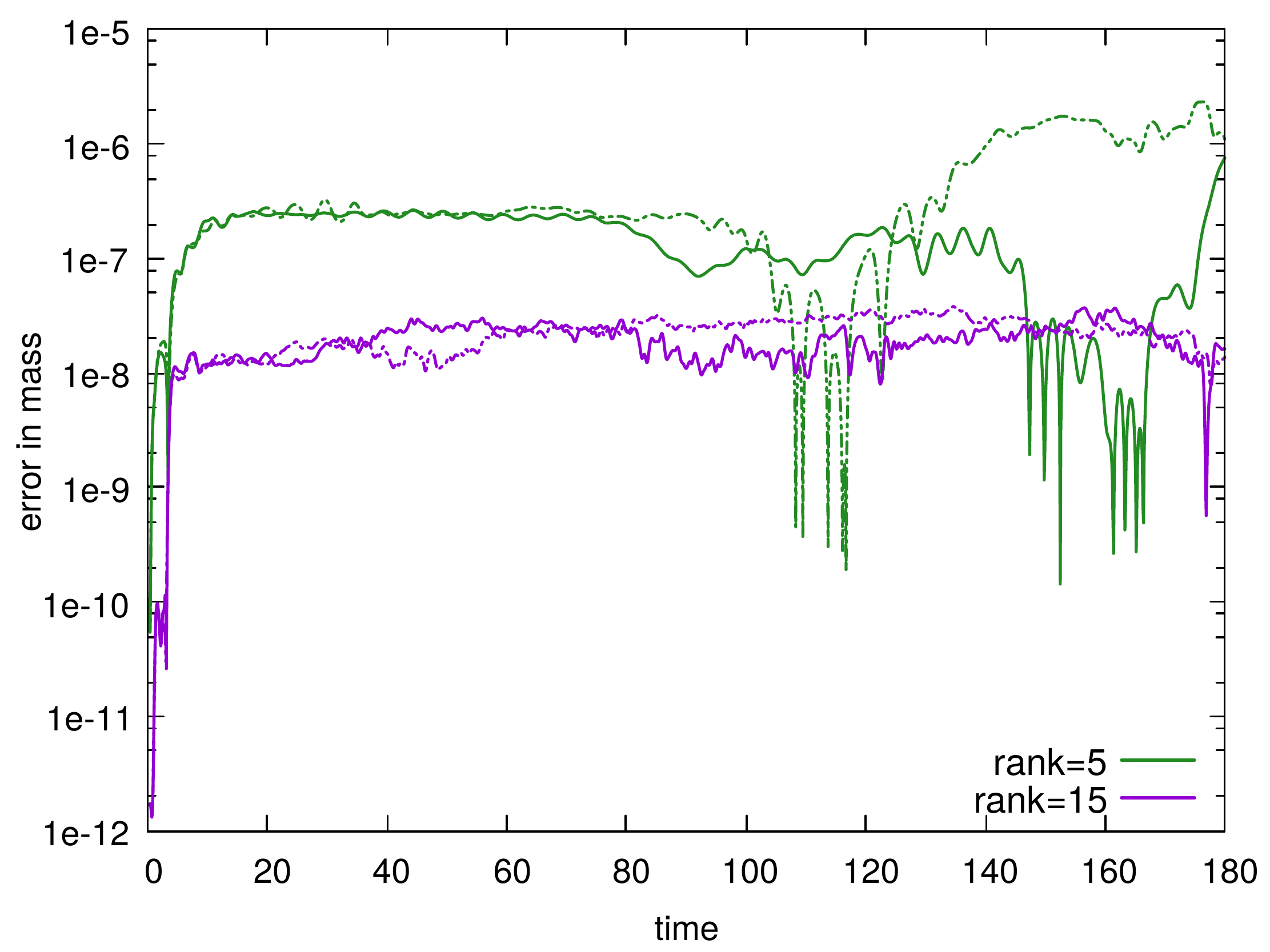}
	\includegraphics[width = .5\textwidth]{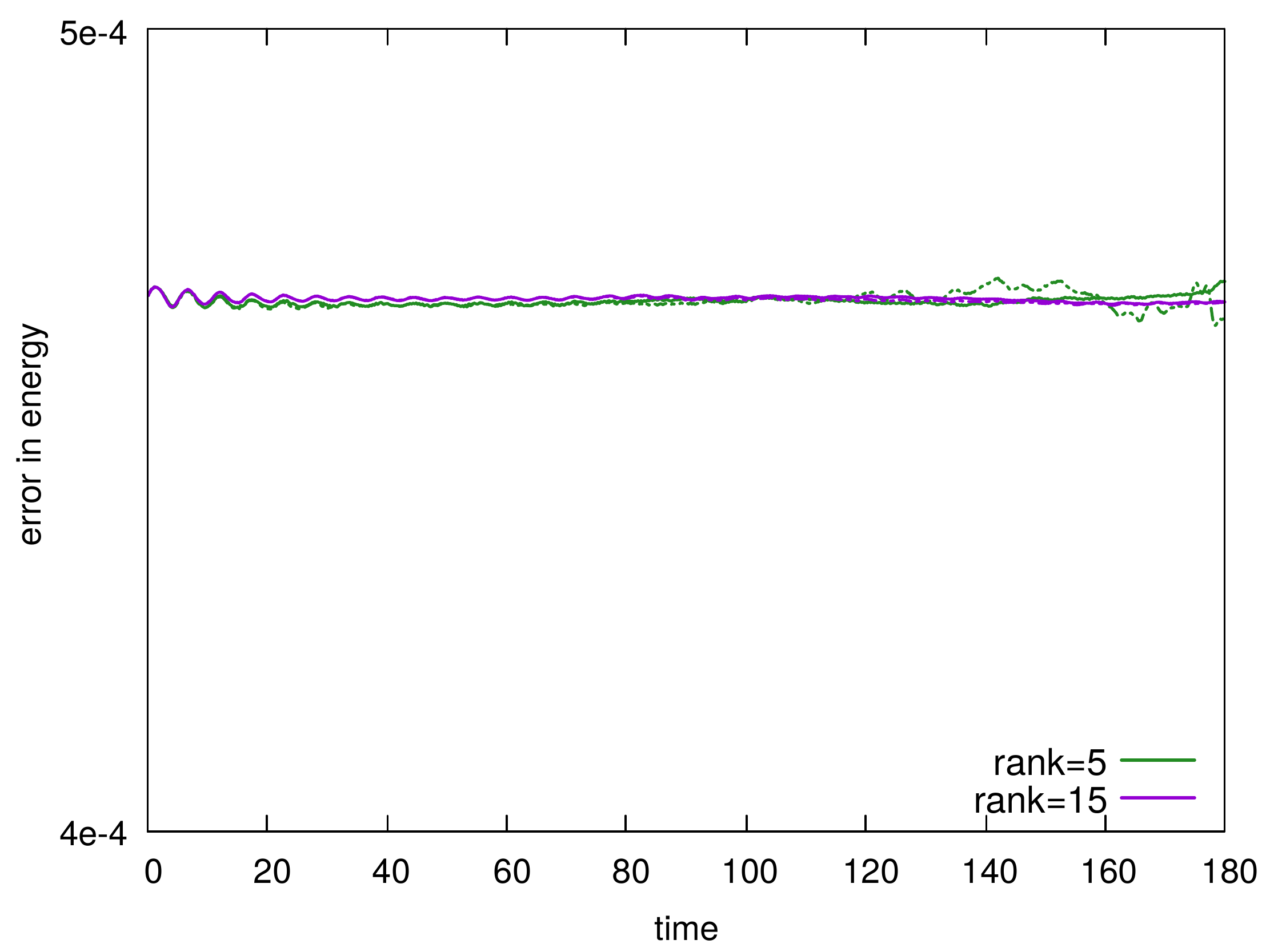} \\
	\center{
	\includegraphics[width = .5\textwidth]{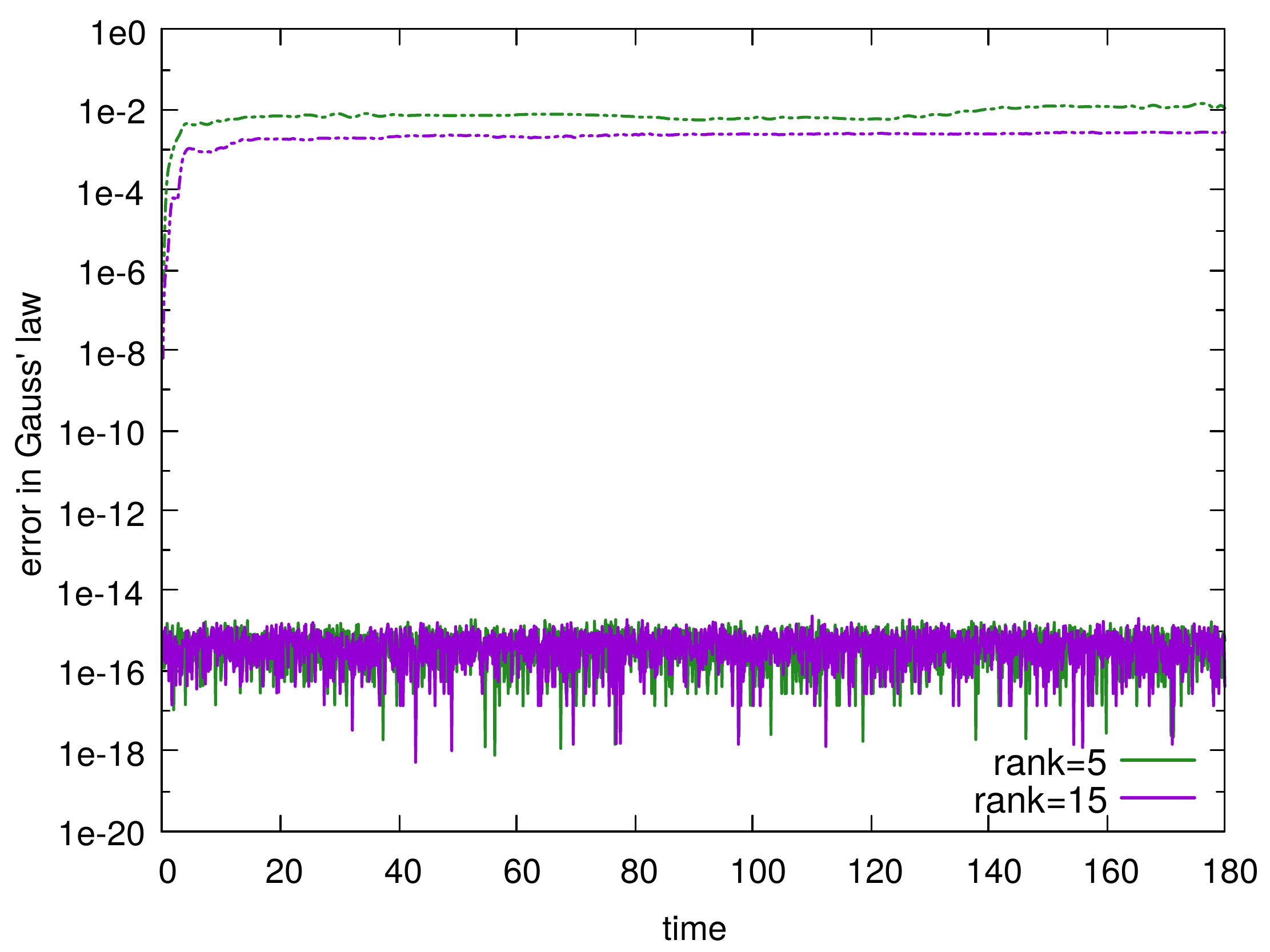}}
	\caption{Landau-type problem. Results for $n_{x_1}=33$, $n_{v_1}=128$, $n_{v_2}=128$, and $\tau=0.1$ (semi-log scale). We compare the low-rank integrator described in algorithm \ref{alg:Strang} (without correction, dashed-dotted lines) with the divergence preserving scheme given in algorithm \ref{alg:Strang_corr} (full lines);  top left: error in mass as defined in \eqref{eq:err_mass}, top right: error in energy as defined in \eqref{eq:err_energy}, bottom: error in Gauss' law measured measured in the discrete $L^2$ norm.}
	\label{fig:lt_comp}
\end{figure}

In figure \ref{fig:lt_comp} we compare the behaviour of the integrator given in algorithm \ref{alg:Strang} (dashed-dotted lines) and the divergence preserving scheme of algorithm \ref{alg:Strang_corr} (full lines) in terms of mass and energy conservation, and preservation of Gauss' law.
The first plot shows the relative error in mass 
\begin{equation} \label{eq:err_mass}
\lvert m(t) - m(0) \rvert / m(0),
\end{equation}
where the mass $m$ is defined as  
\[ 
m(t) = \int_{\Omega_{x_1}} \int_{\Omega_{v}} f(t,x_1,v) \dd x_1 \dd v. 
\] 
The relative error in the total energy $\mathcal{E}_ {tot} =  \mathcal{E}_{e}(t) + \mathcal{E}_{m}(t) + \mathcal{E}_{k}(t)$, i.e. 
\begin{equation} \label{eq:err_energy}
\lvert \mathcal{E}_{tot} (t) - \mathcal{E}_{tot}(0) \rvert / \mathcal{E}_{tot}(0),
\end{equation}
is shown in the second plot. 
Conservation of mass and energy is quite good for both algorithms. For the case of algorithm \ref{alg:Strang} we observe that mass conservation and the error in Gauss' law improve  as we increase the rank indicating that the error is dominated by the low-rank approximation. We also observe that among the three invariants considered here, the largest error, approximately $10^{-3}$ for rank $15$, is committed in Gauss' law. When applying the divergence preserving scheme of algorithm \ref{alg:Strang_corr} Gauss' law is satisfied up to machine precision. We further note that the divergence correction does not cause a significant difference in the behavior of the other invariants.

\section{Numerical results for the two-stream instability} \label{sect:two-stream}
In this section we consider the two-stream instability. We employ the problem found in \cite{CEF15}. 
The initial particle density is given by 
\[ 
f(0,x_1,v_1,v_2) = \frac{1}{2\pi \beta} \ee^{-v_2^2/\beta}\left(\ee^{-(v_1-0.2)^2/\beta}+\ee^{-(v_1+0.2)^2/\beta}\right),
\]
where we choose $\beta = 2 \cdot 10^{-3}$, $\Omega_{x_1} = [0,2\pi]$ and $\Omega_{v_1}\times \Omega_{v_2} = [-0.4,0.4]^2$. The instability is driven by a magnetic perturbation only. The initial magnetic field is chosen as 
\[ 
B_3(0,x_1) = \alpha \sin x_1
\]
with $\alpha = 10^{-3}$. The electric field is initialized to zero.

\begin{figure}[h]
	\includegraphics[width = .5\textwidth]{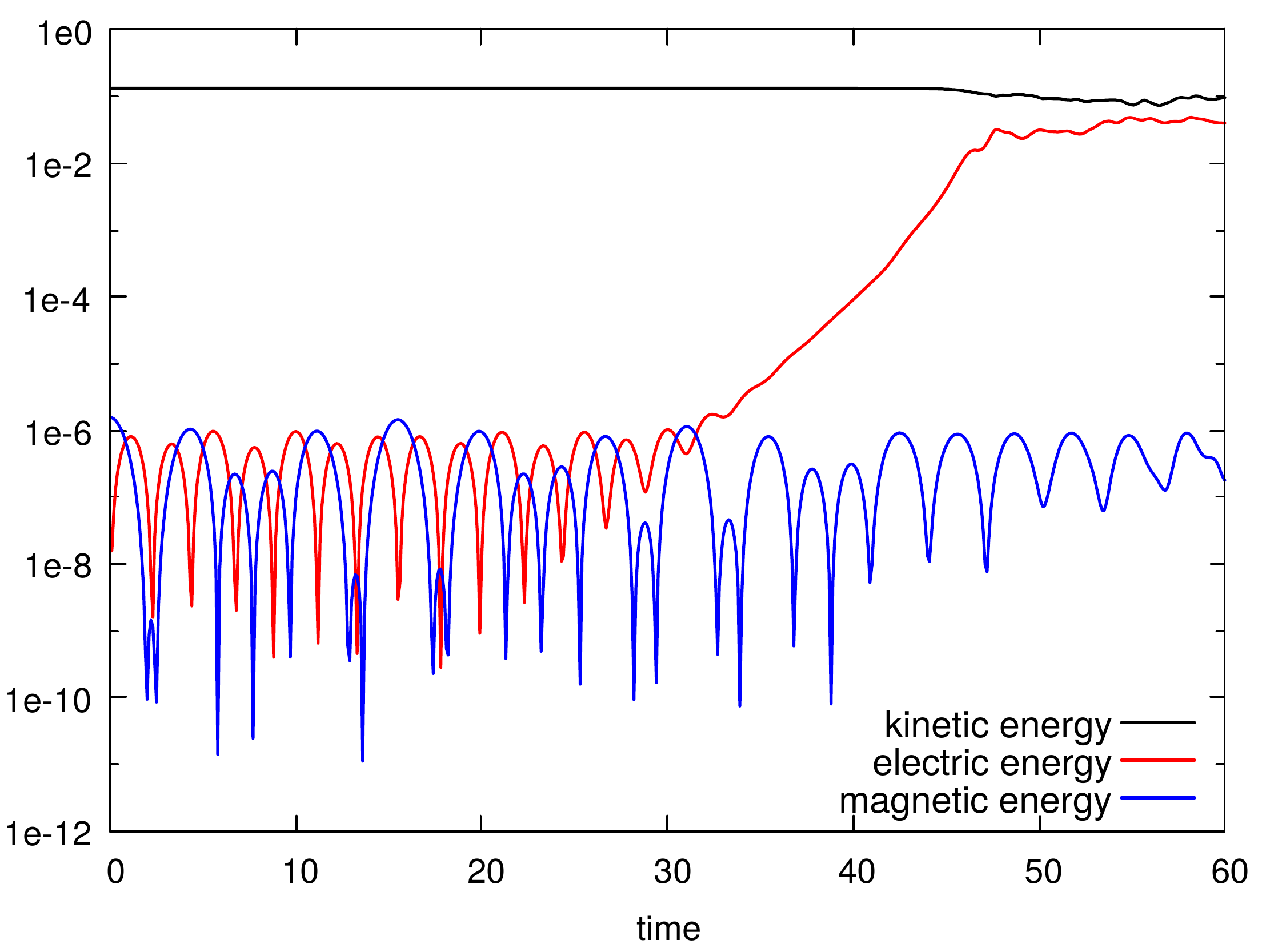}
	\includegraphics[width = .5\textwidth]{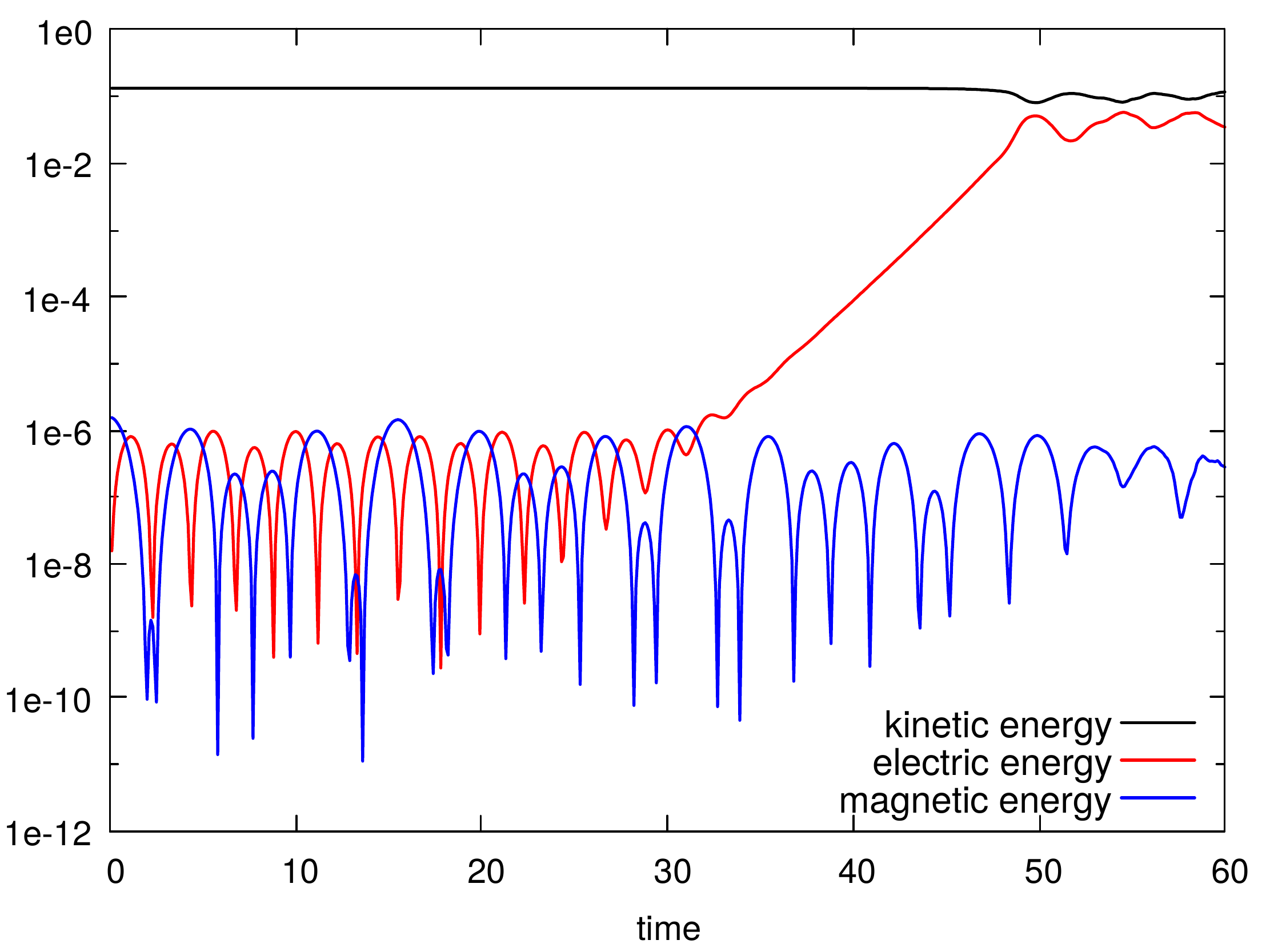}
	\includegraphics[width = .5\textwidth]{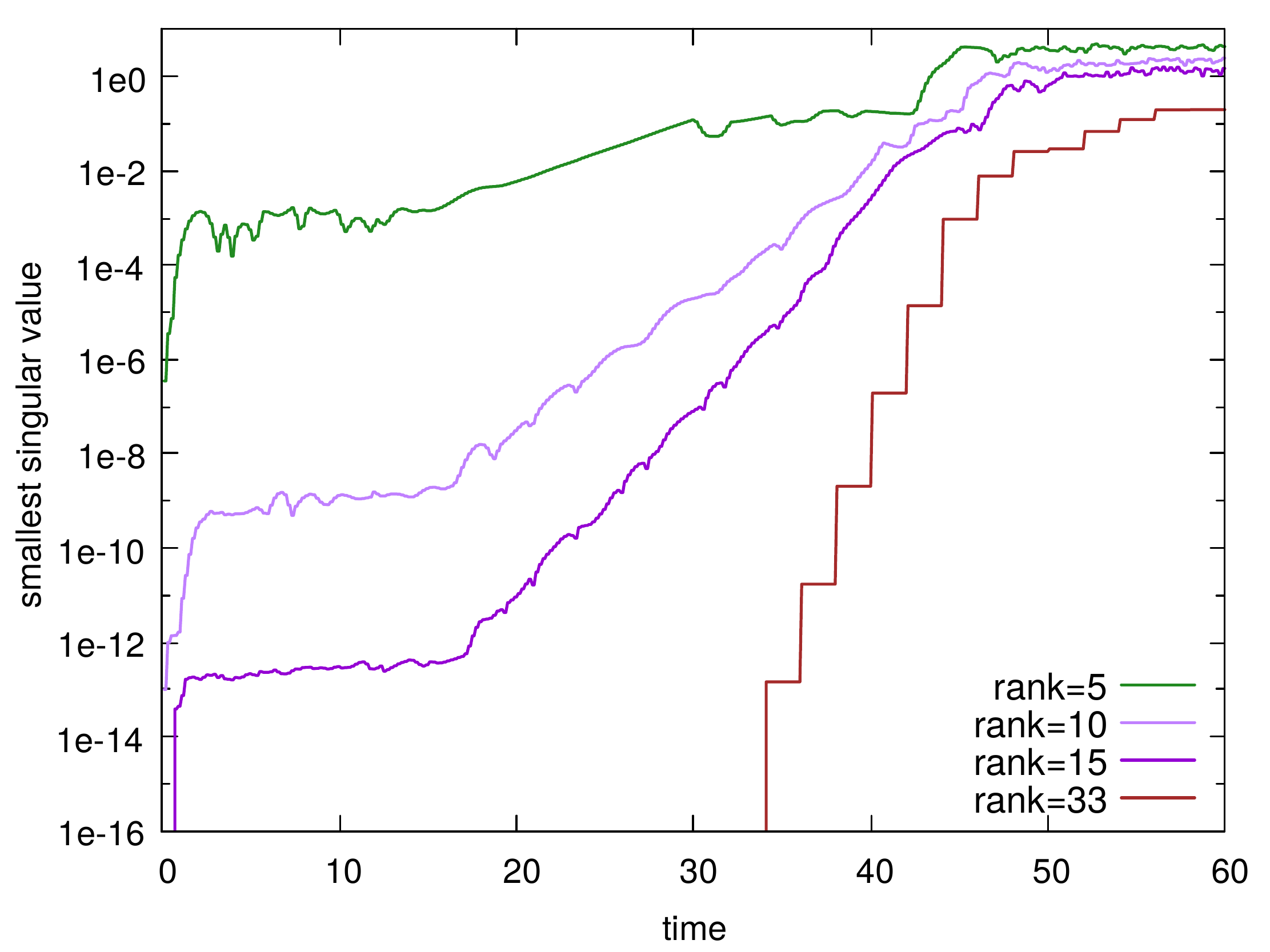}
	\includegraphics[width = .5\textwidth]{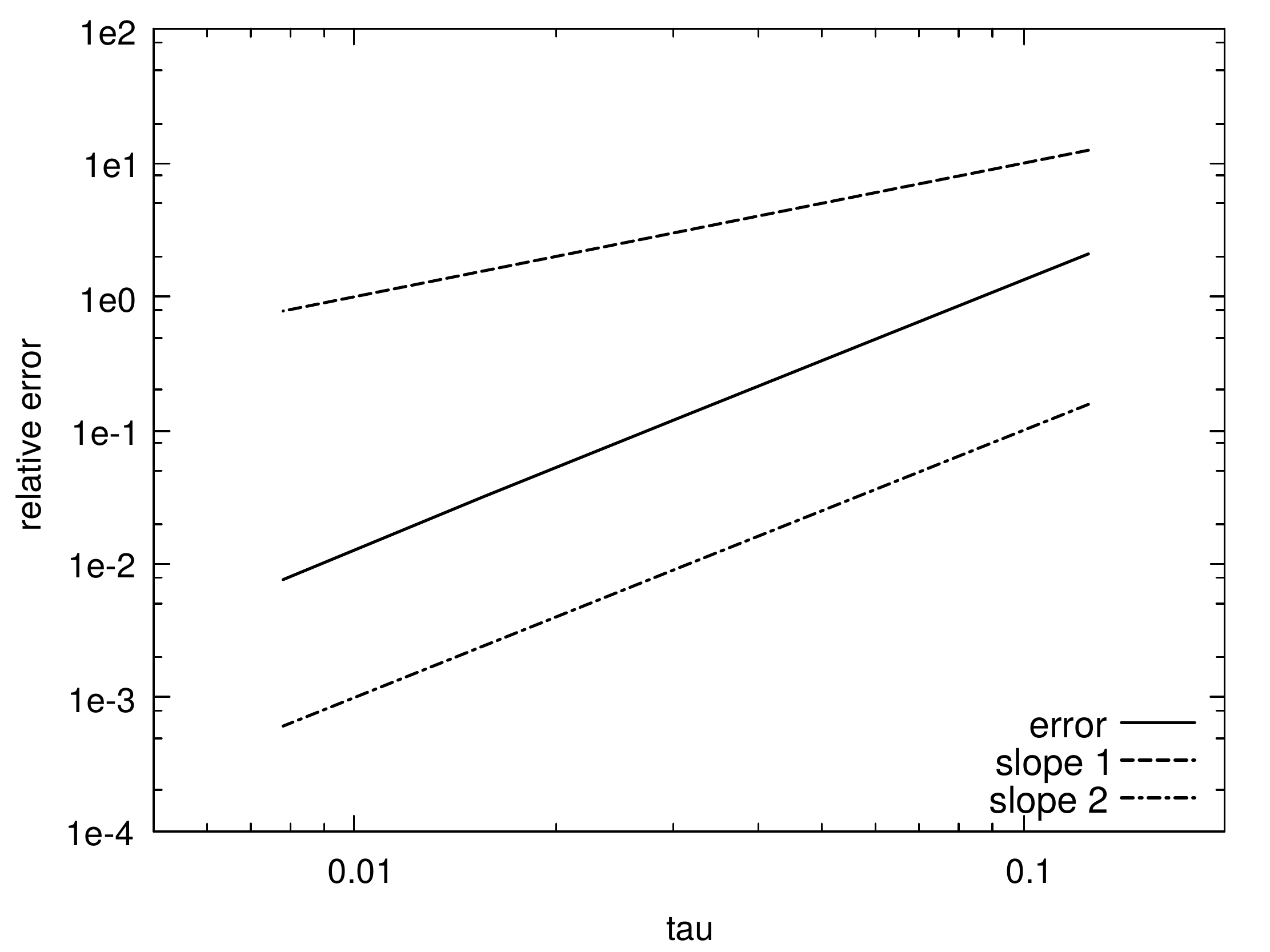} 
	\caption{Two-stream instability. Results for algorithm \ref{alg:Strang} (without correction) for $n_{x_1}=33$, $n_{v_1}=64$, and $n_{v_2}=64$. Top  left: energies for rank = 5 and $\tau=0.1$ (semi-log scale); top right: energies for rank = 15 and $\tau=0.1$ (semi-log scale); bottom left: smallest singular value of the corresponding distribution function computed every 20 time steps (semi-log scale); bottom right: relative error in the distribution function at $T=2$ for different values of $\tau$ and rank = 15 (log-log scale).} 
	\label{fig:2si_wc_energy}
\end{figure}

In figure \ref{fig:2si_wc_energy}, on the top, we plot the time evolution of the energies obtained with algorithm \ref{alg:Strang} (without correction) for different values of the approximation rank. 
Initially, we observe an oscillatory behavior of the electric and magnetic energies. Later, the electric energy increases exponentially and then saturates. 
The results for rank = 5 are on the left, on the right the results for rank = 15 are displayed.
We clearly notice the benefit of choosing a higher approximation rank. In particular, the electric and kinetic energies show a better match after saturation. This might indicate that a higher approximation rank is necessary to resolve the nonlinear regime. 
On the other hand, the plot of the smallest singular values brings a different conclusion. The error in the distribution function in the nonlinear regime is relatively large (independent of the choice of the rank). Despite this large error in the distribution function, we observe that averaged quantities of the system, such as the electric energy, match the physical dynamics quite well.

In figure \ref{fig:2si_wc_energy} we also show that the splitting integrator described in algorithm \ref{alg:Strang} converges in time with the expected order. We display the time error of the scheme for different values of the time step $\tau$ with respect to a reference solution.  

\begin{figure}
	\includegraphics[width = .5\textwidth]{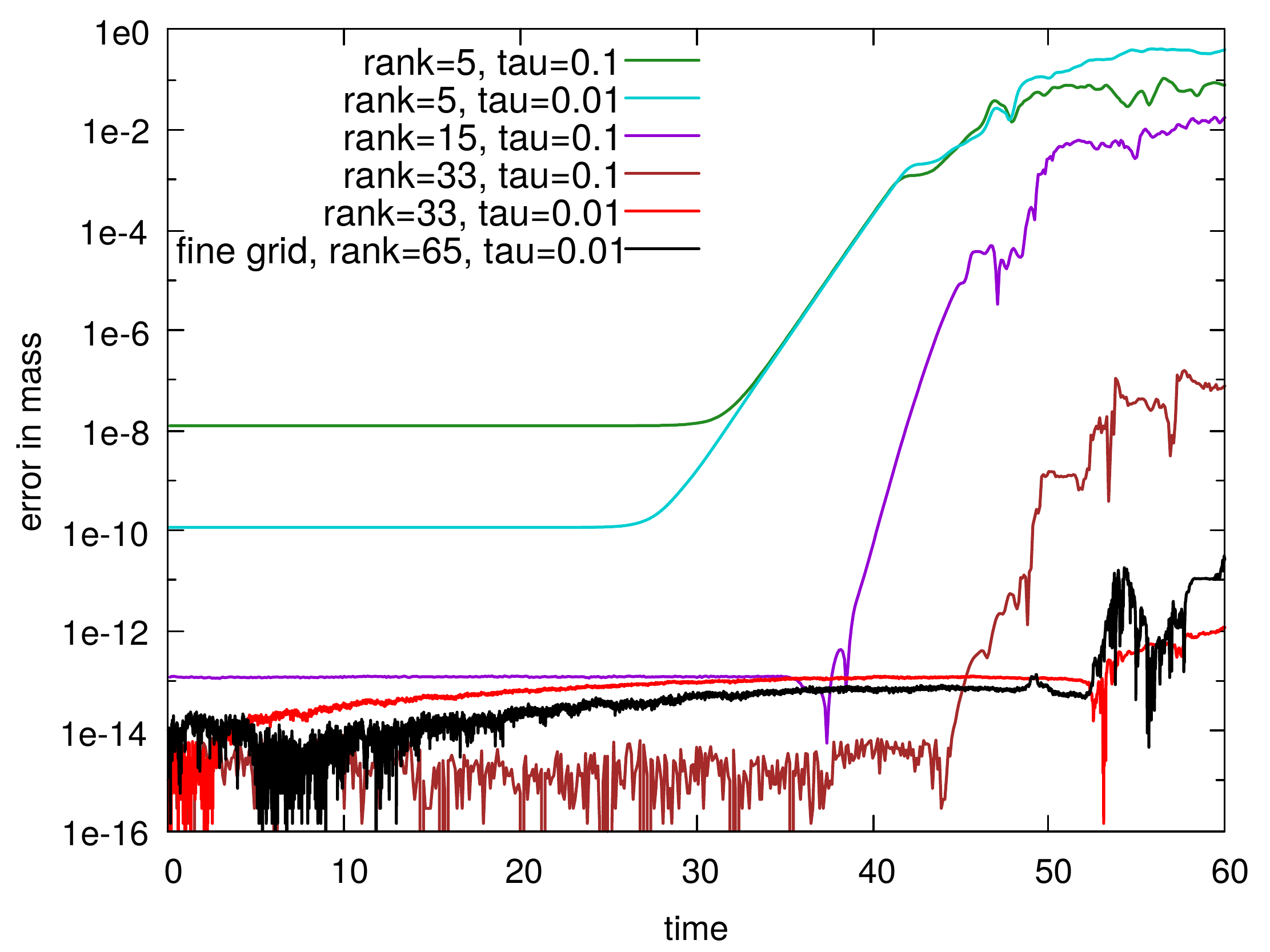}
	\includegraphics[width = .5\textwidth]{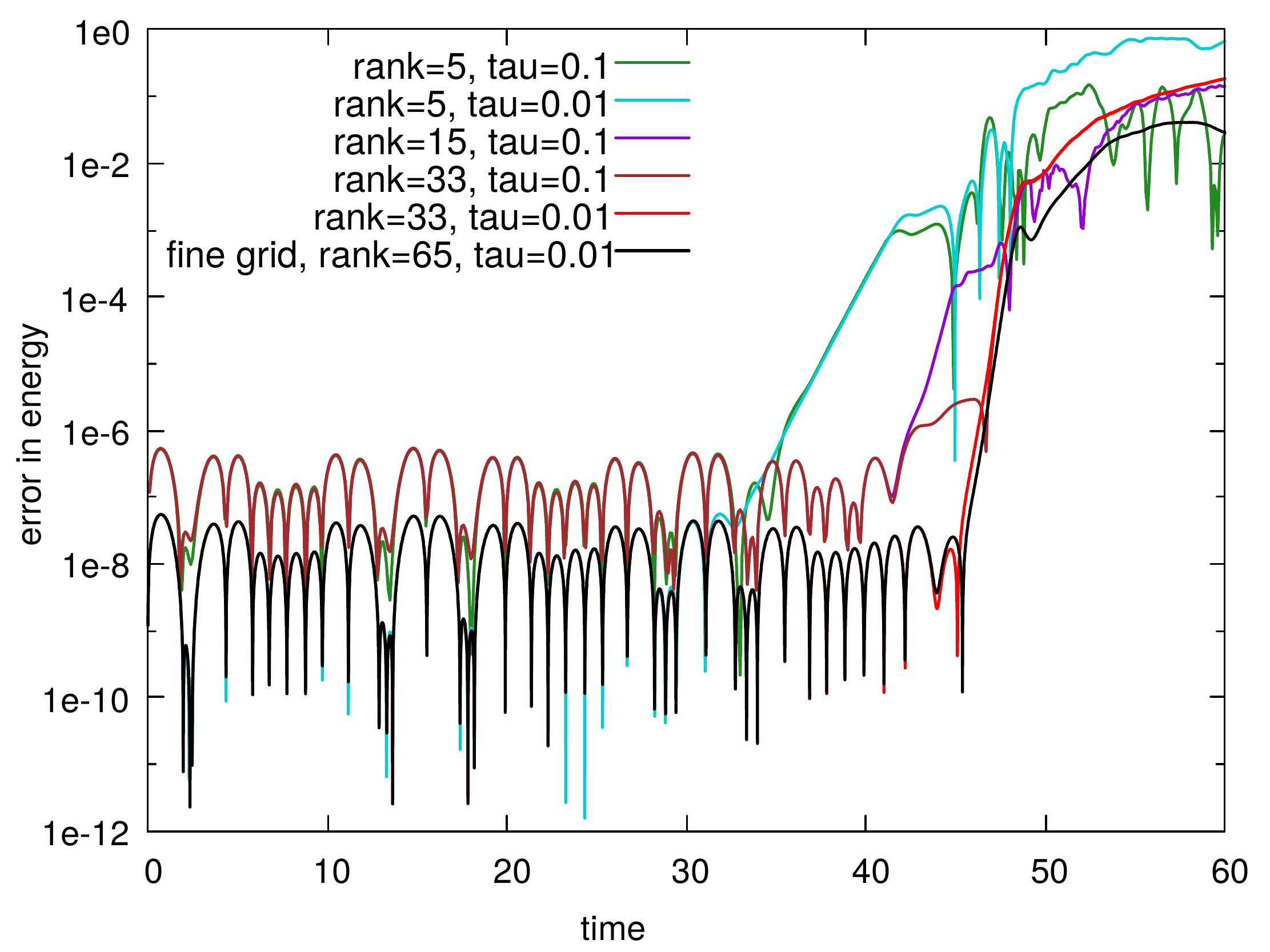} \\
	\center{
	\includegraphics[width = .5\textwidth]{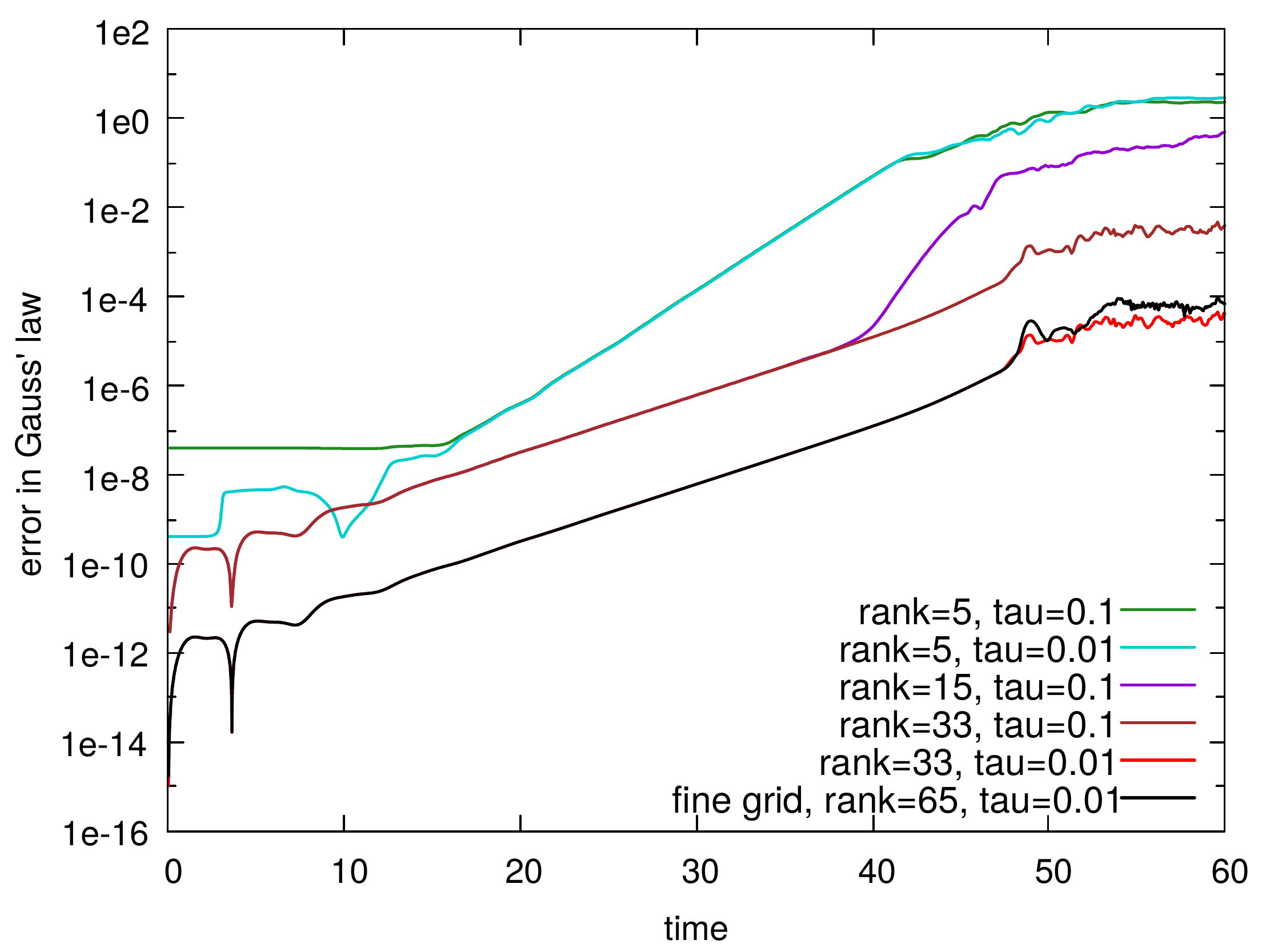}}
	\caption{Two-stream instability. Results for algorithm \ref{alg:Strang} (without correction) (semi-log scale). Unless otherwise specified the space discretization is $n_{x_1}=33$, $n_{v_1}=64$, and $n_{v_2}=64$. The fine grid is obtained with $n_{x_1}=65$, $n_{v_1}=128$, and $n_{v_2}=128$.  top left: error in mass as defined in \eqref{eq:err_mass}; top right: error in energy as defined in \eqref{eq:err_energy}; bottom: error in Gauss' law measured measured in the discrete $L^2$ norm.}
	\label{fig:2si_wc_comp}
\end{figure}

In figure \ref{fig:2si_wc_comp} we show the behavior of the integrator with respect to the conservation of mass, energy, and the preservation of Gauss' law.  The numerical error of the proposed integrator is an interplay between different factors. We consider the influence of the choice of the approximation rank, the spatial discretization, and the time integration error due to the projector splitting integrator. 
On the left, the error in mass is plotted. First, let us focus on the curves corresponding to different approximation ranks obtained with time step size $\tau = 0.1$. We observe that increasing the rank gives more accurate results. In particular, the method preserves the mass up to machine precision for the full rank (i.e. rank = 33 for the configurations with $n_{x_1}=33$) until saturation. Note, however, that in the nonlinear regime the error increases in time. This is due to the fact that even if the dynamics can be exactly represented by the chosen rank, there is en error comitted by the projector splitting integrator. In this regime, conservation of mass can be enforced by choosing a smaller time step (we consider $\tau = 0.01$ in the figure).

Similar conclusion can be drawn by looking at the error in energy shown on the right. In this case we can further observe that a finer space discretization is beneficial in the nonlinear regime. In fact, the red and black curves overlap until saturation. Afterwards, the error for the fine grid is smaller compared to the other configurations. The plot of the error in Gauss' law, also, confirms the previous observations. If the approximation rank is high enough, the time integration error can be clearly observed. This error decreases by reducing the time step size from $\tau =0.1$ and $\tau = 0.01$ by two orders of magnitude, as expected for a second-order integrator. 

\begin{figure}
	\includegraphics[width = .5\textwidth]{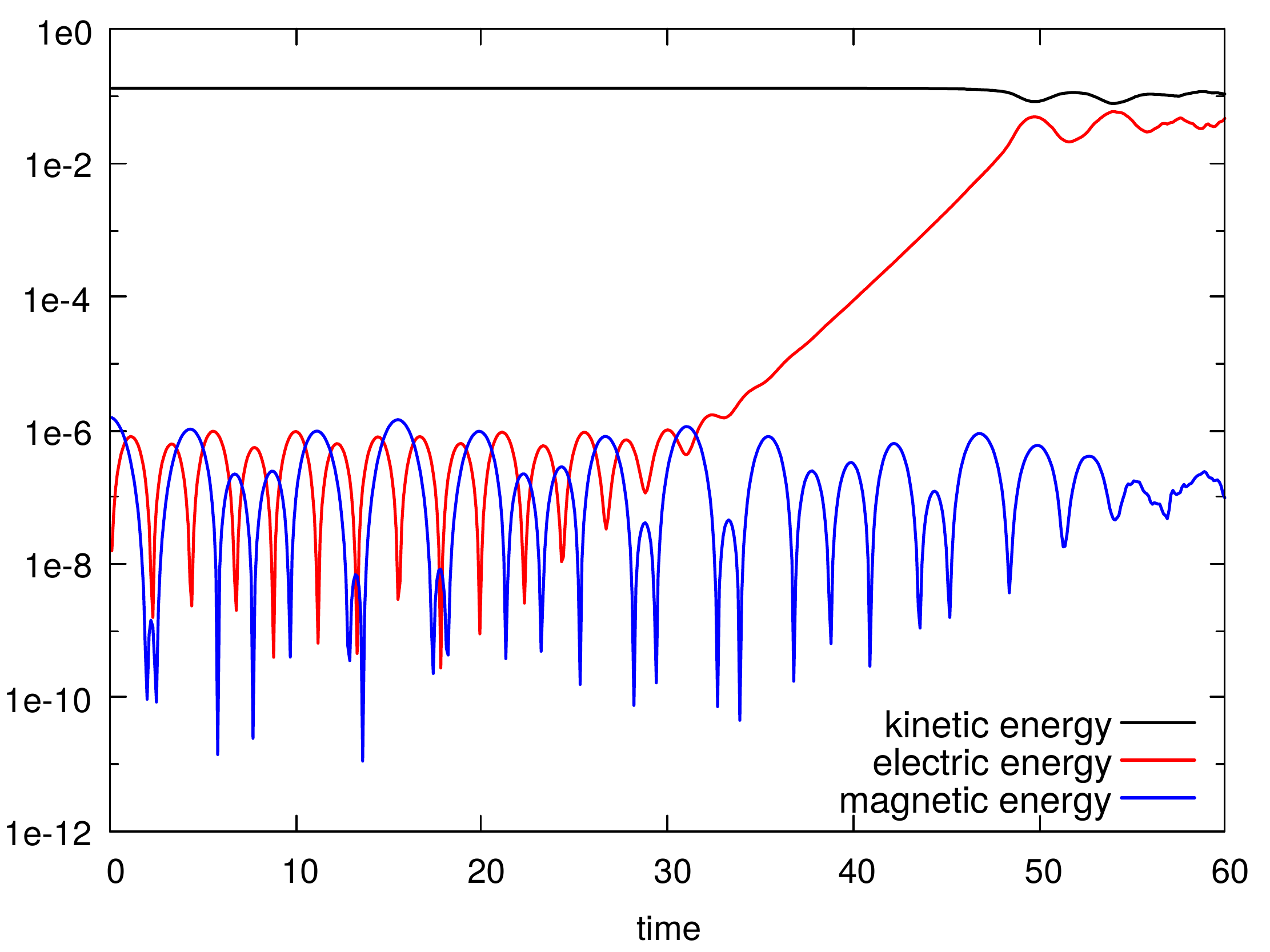} 
	\includegraphics[width = .5\textwidth]{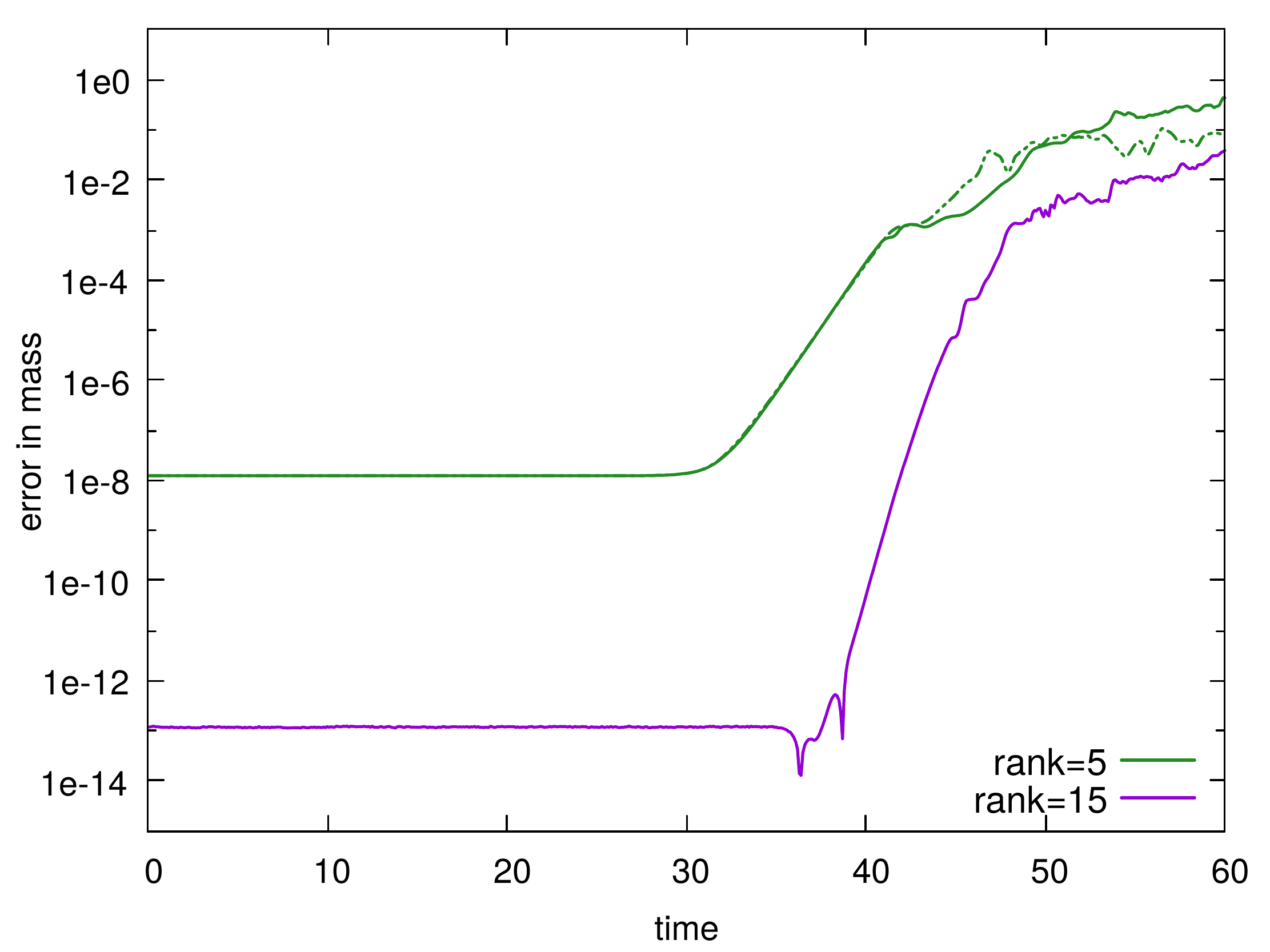}
	\includegraphics[width = .5\textwidth]{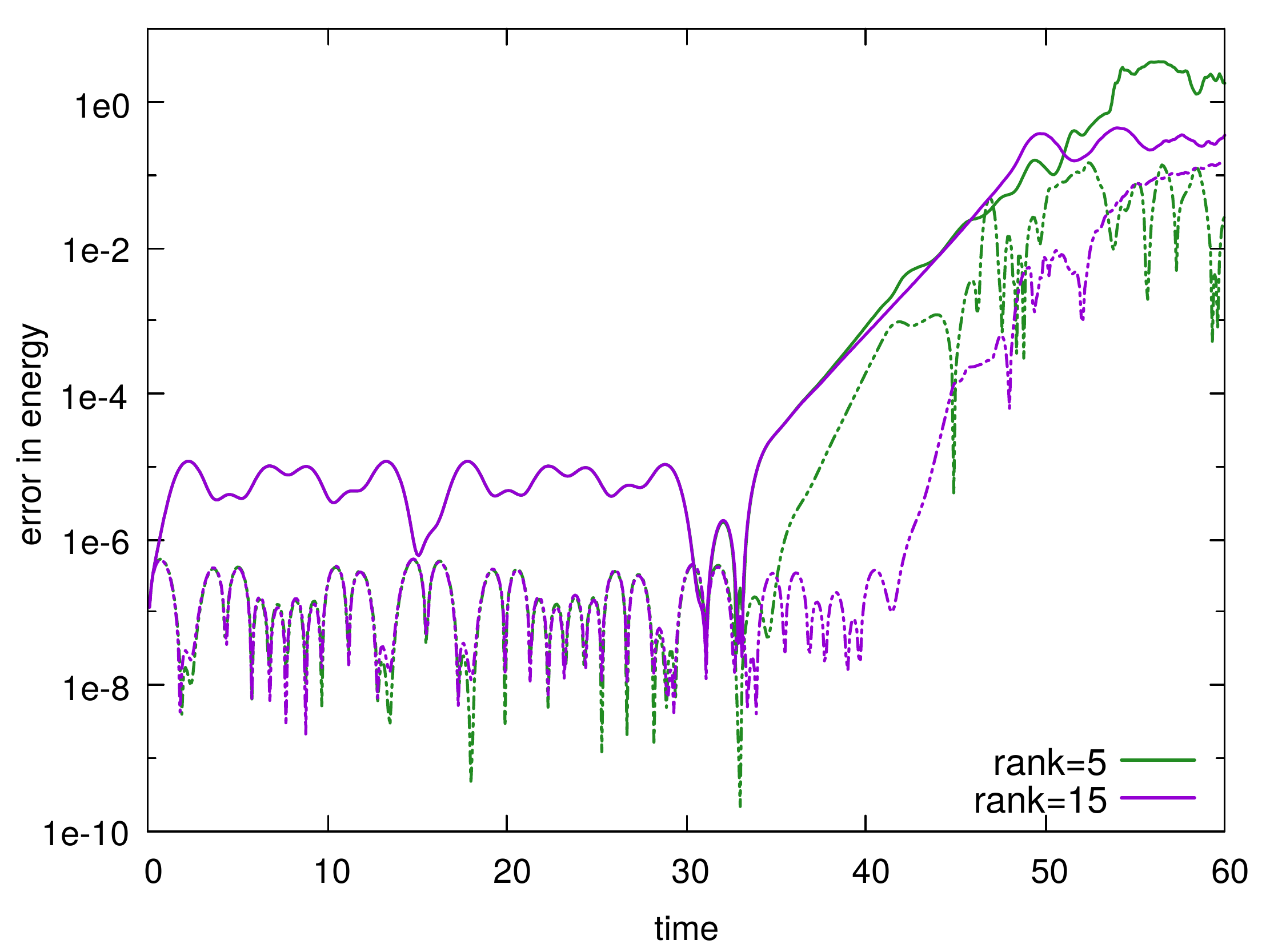}
	\includegraphics[width = .5\textwidth]{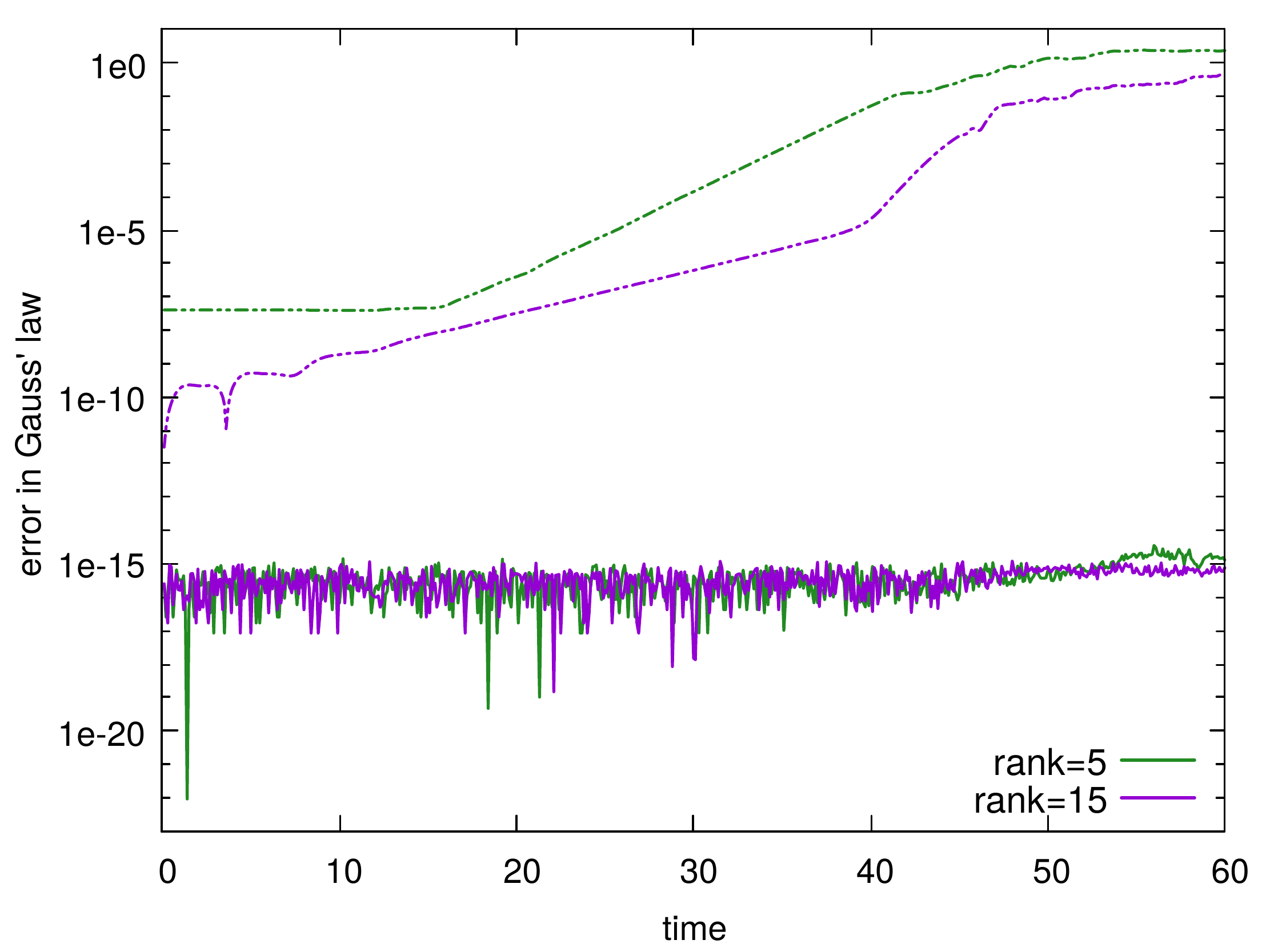}
	\caption{Two-stream instability. Results for $n_{x_1}=33$, $n_{v_1}=64$, $n_{v_2}=64$, and $\tau = 0.1$ (semi-log scale). Top left: energies for rank = 15 for algorithm \ref{alg:Strang_corr} (with correction). In the other figures we compare the low-rank integrator described in algorithm \ref{alg:Strang} (without correction, dashed-dotted lines) with the divergence preserving scheme given in algorithm \ref{alg:Strang_corr} (full lines); top right: error in mass as defined in \eqref{eq:err_mass}, bottom left: error in energy as defined in \eqref{eq:err_energy}, bottom right: error in Gauss' law measured measured in the discrete $L^2$ norm.}
	\label{fig:2si_comparisons}
\end{figure}

In figure \ref{fig:2si_comparisons}, top left, we observe that the time evolution of the energies for the divergence preserving scheme are close to what is obtained with algorithm \ref{alg:Strang}. As for the Landau-type problem, the divergence correction does not cause a significant difference in the behavior of the other invariants.


\section{Numerical results for the bump-on-tail instability} \label{sect:bump}
\begin{figure}[h]
	\subfigure{\includegraphics[height = 0.18\textwidth]{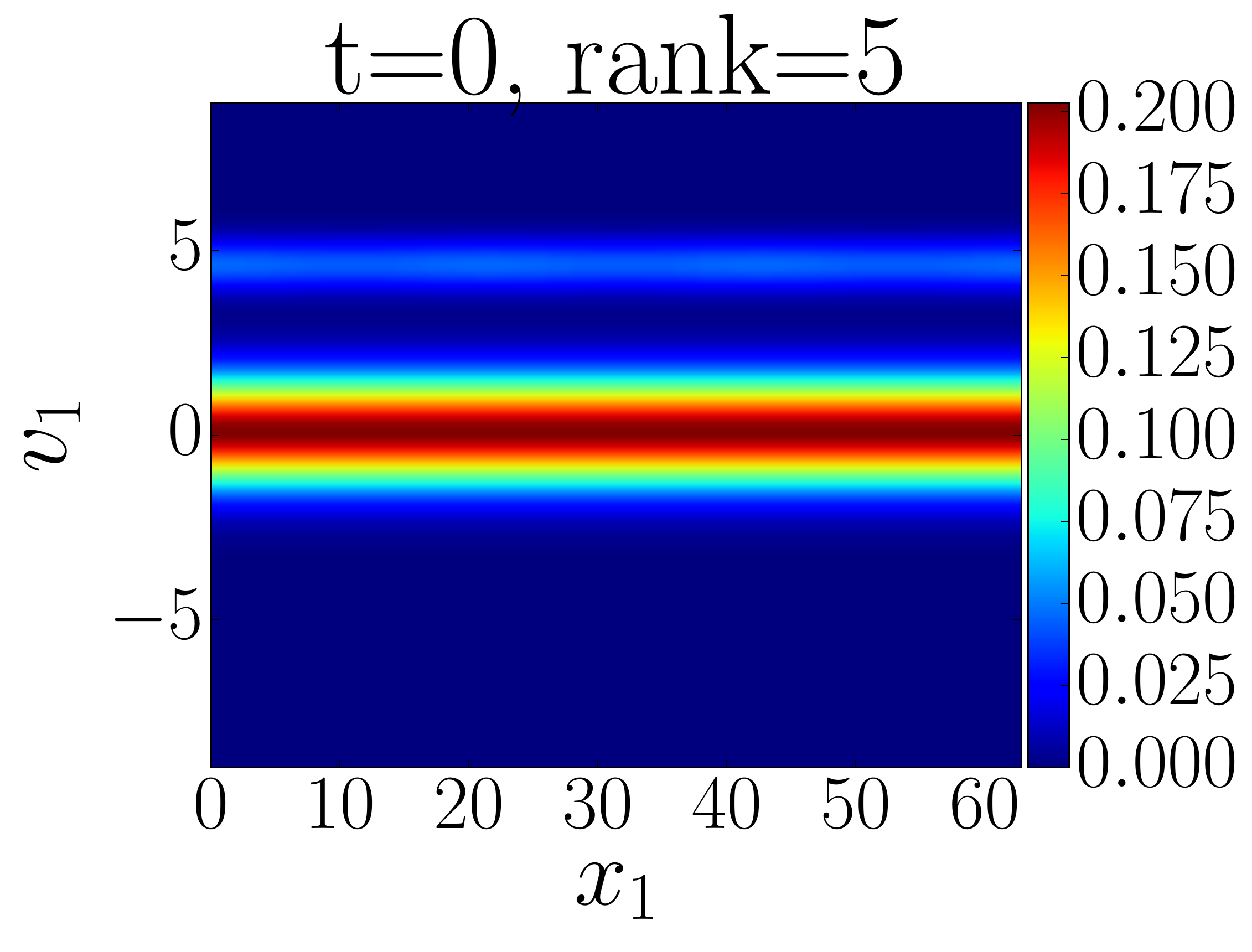}}
	\subfigure{\includegraphics[height = 0.18\textwidth]{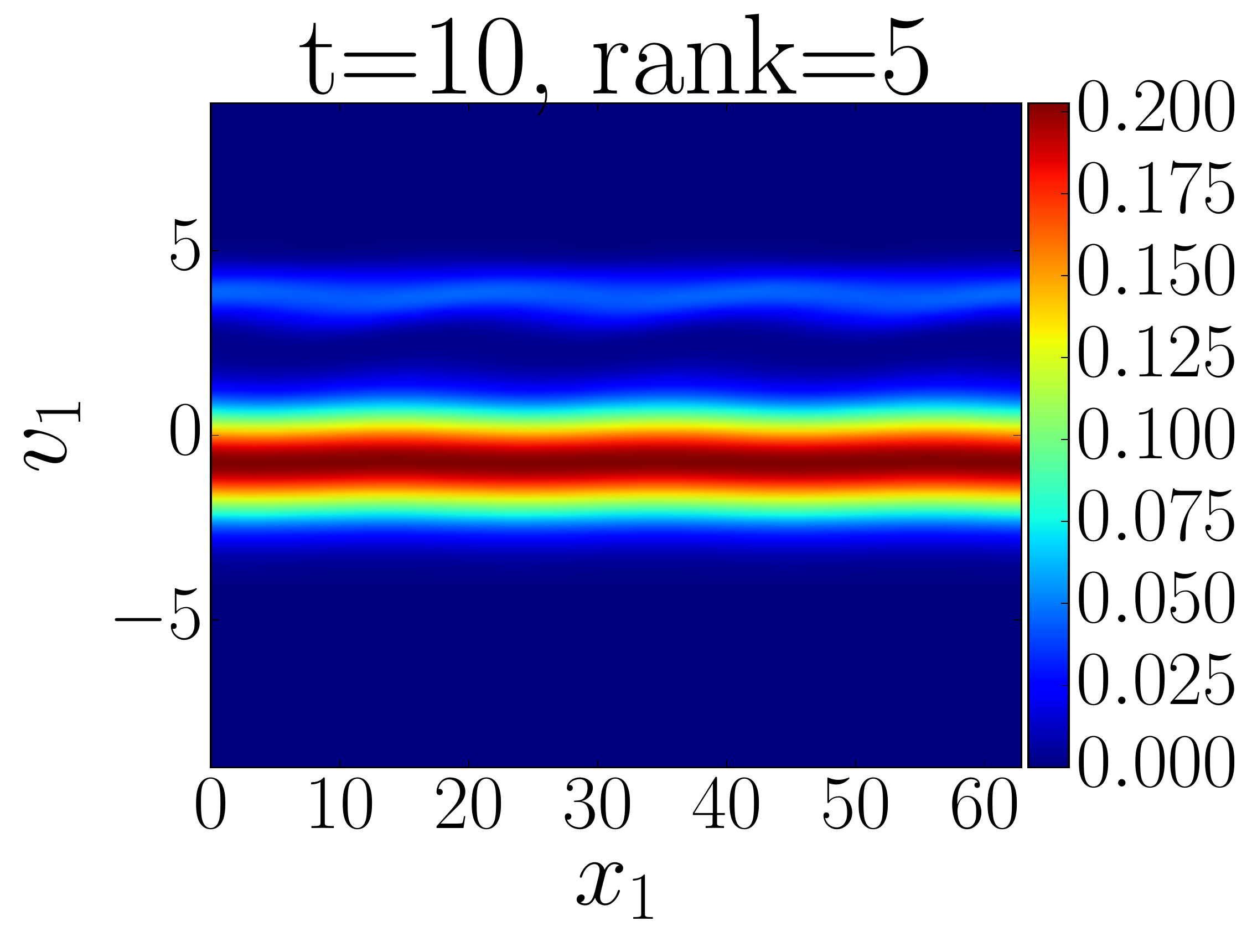}}
	\subfigure{\includegraphics[height = 0.18\textwidth]{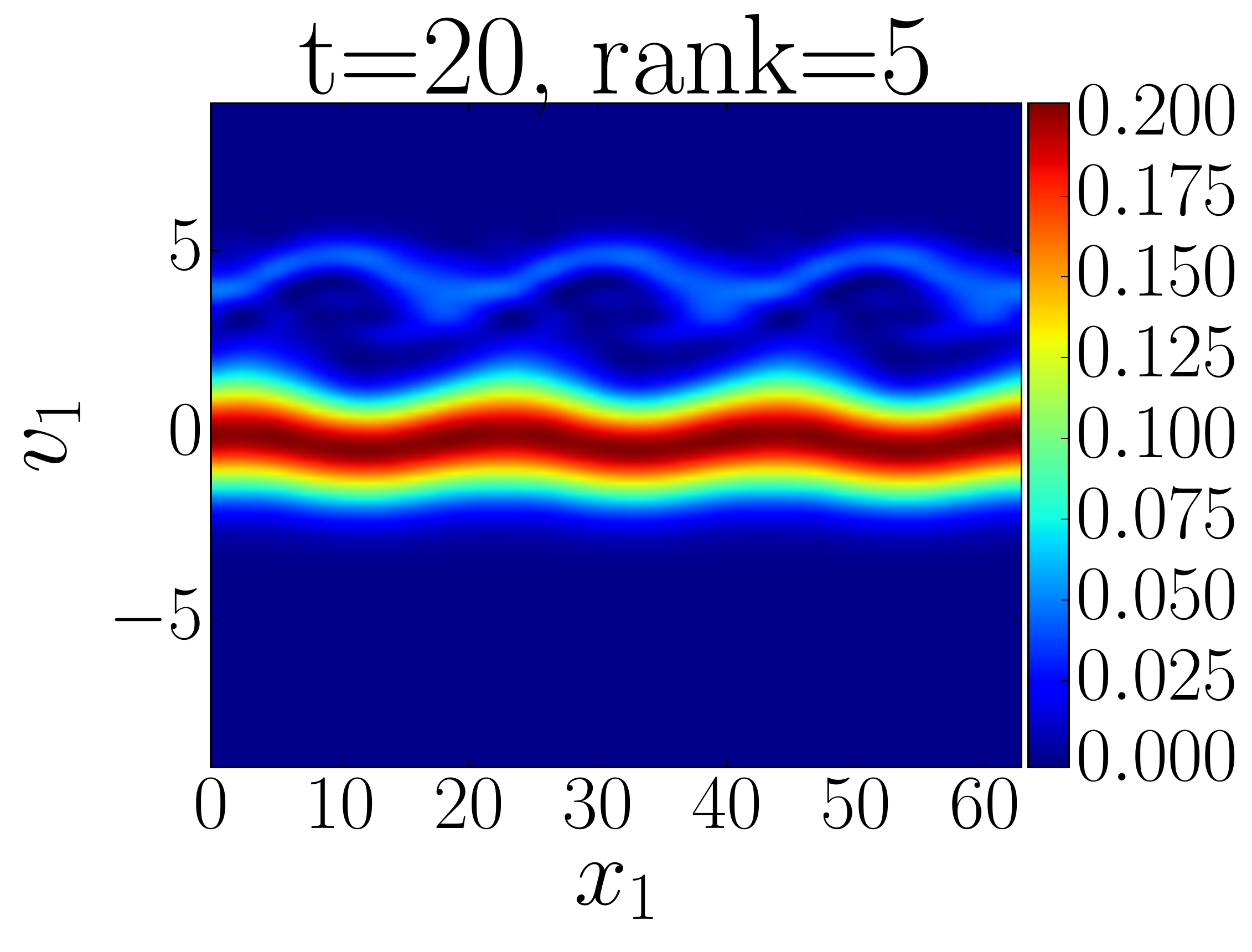}}
	\subfigure{\includegraphics[height = 0.18\textwidth]{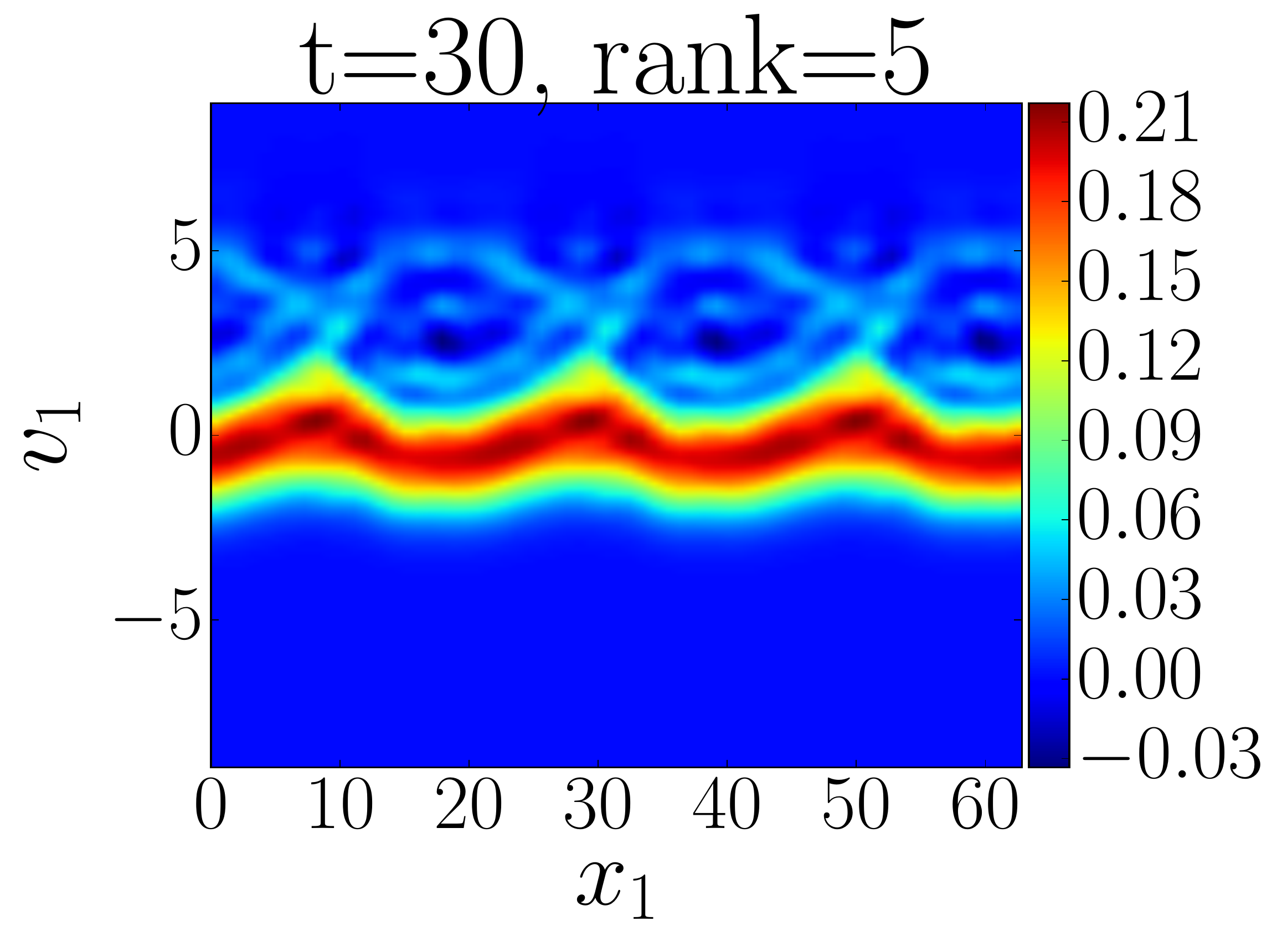}}\\
	\subfigure{\includegraphics[height = 0.18\textwidth]{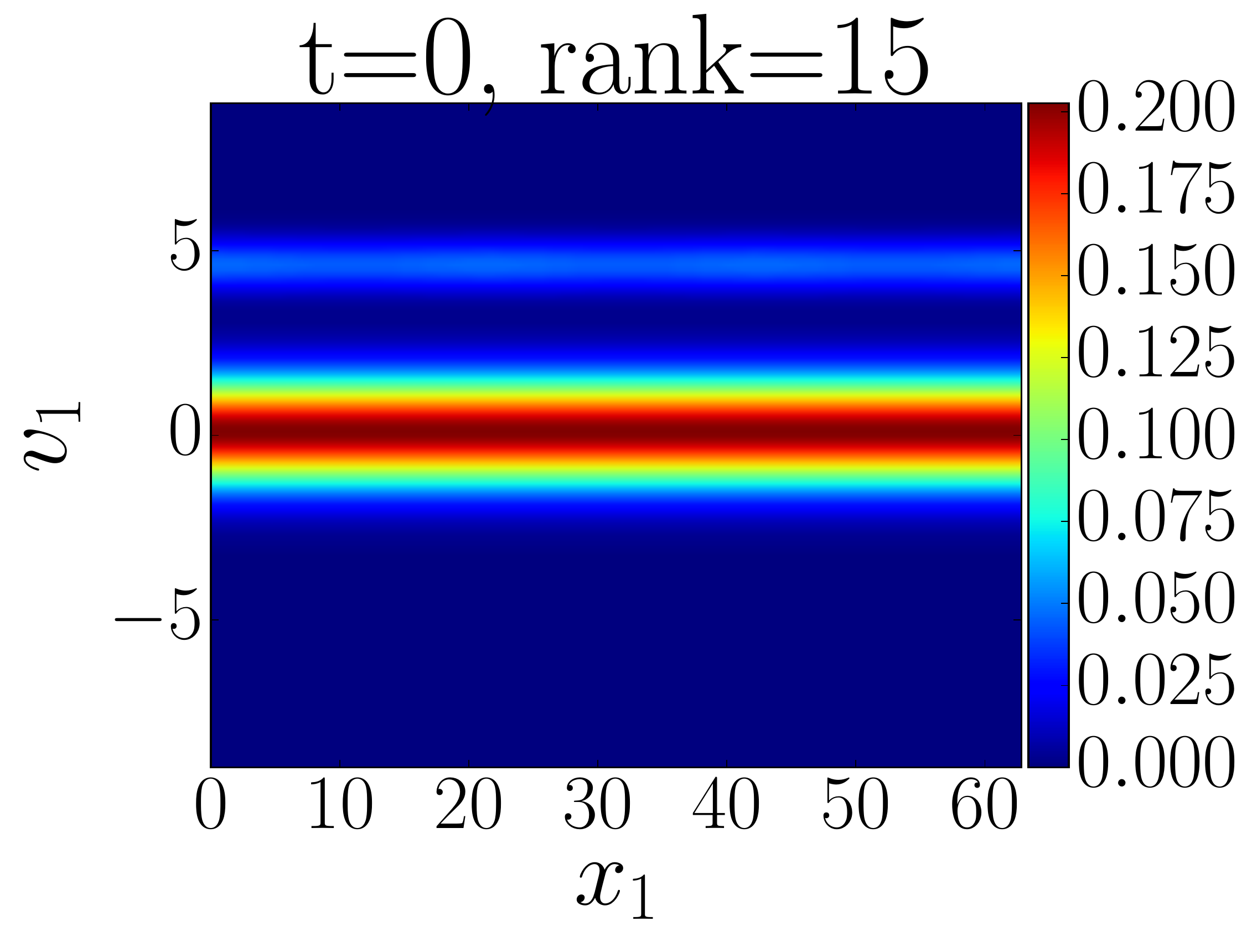}}
	\subfigure{\includegraphics[height = 0.18\textwidth]{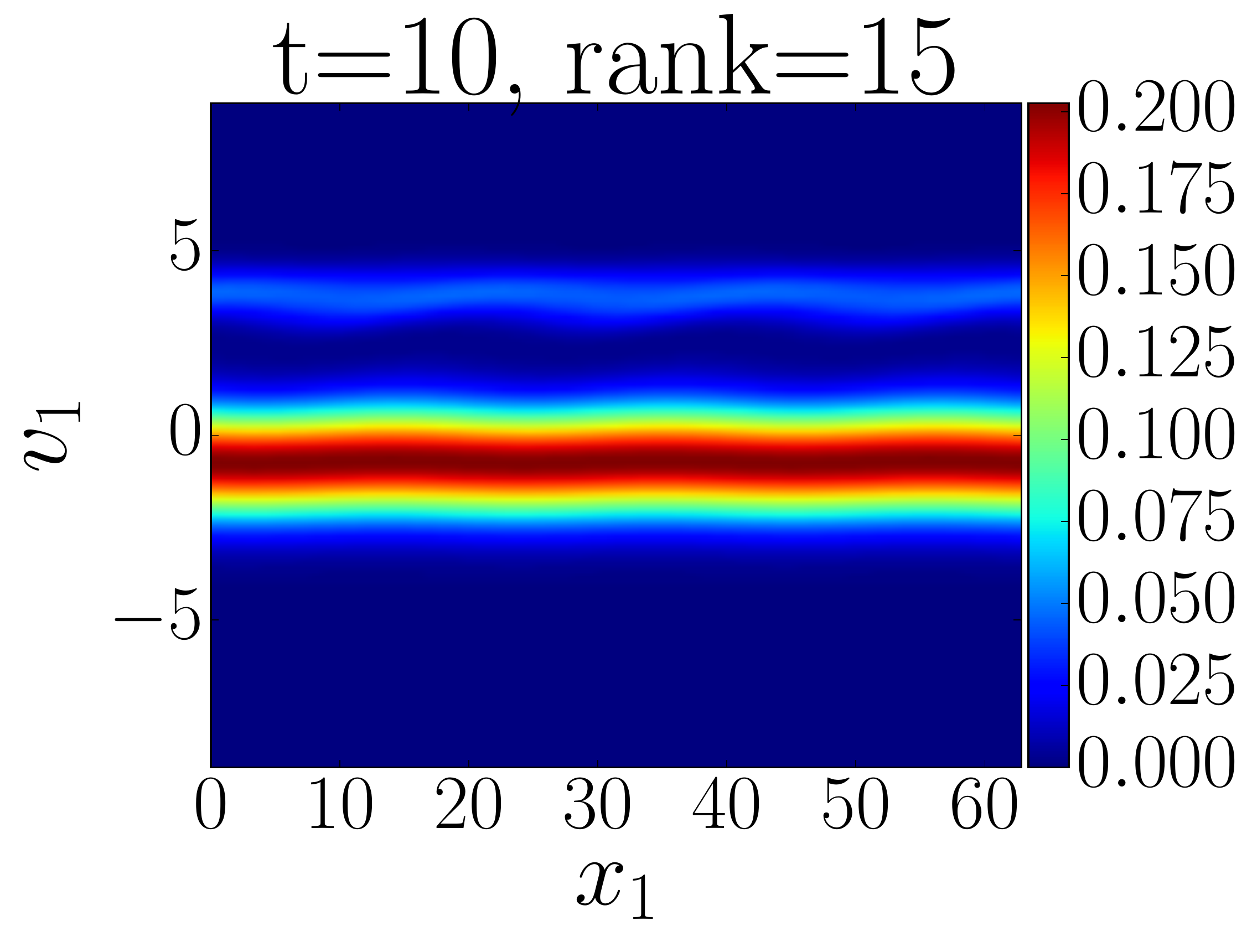}}
	\subfigure{\includegraphics[height = 0.18\textwidth]{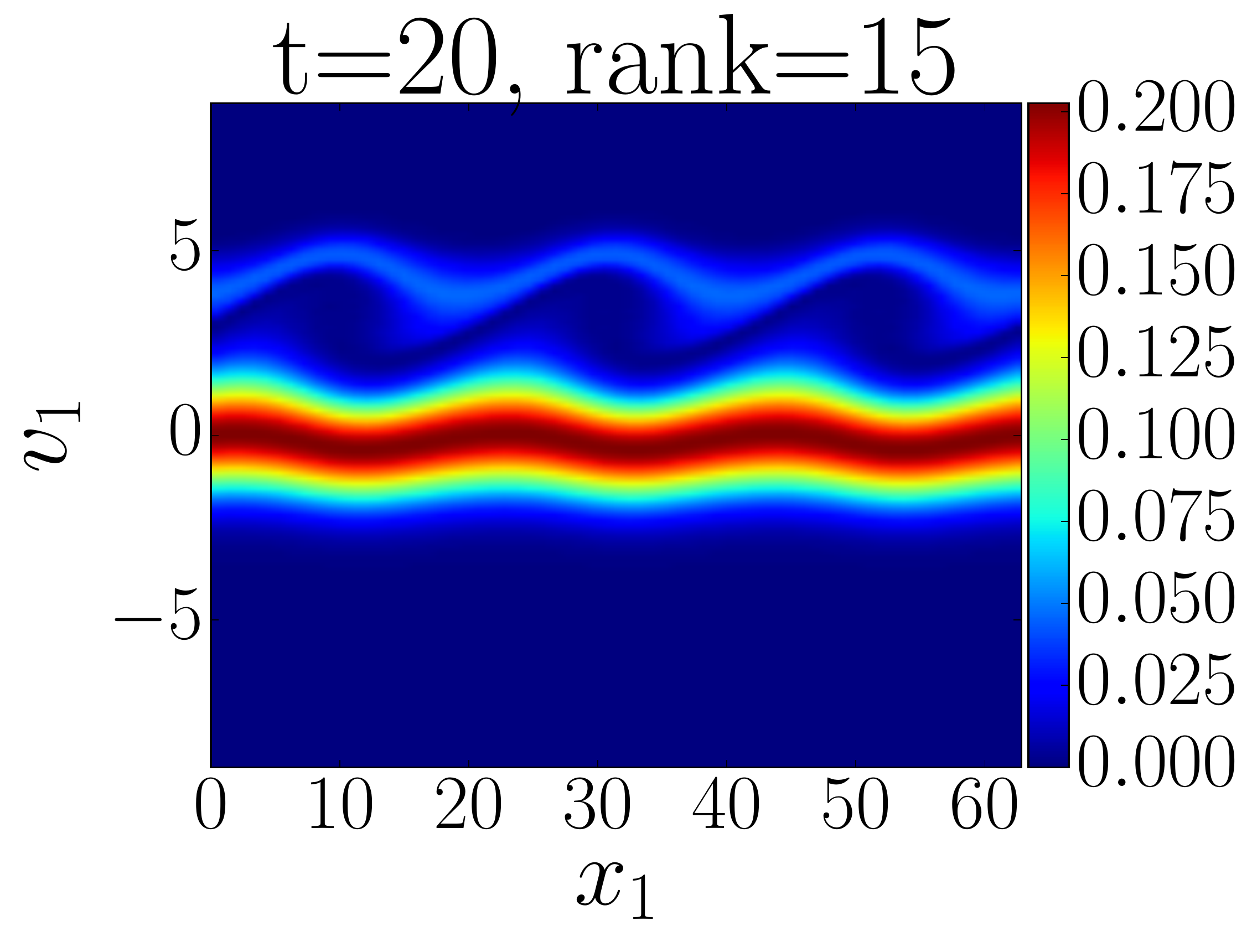}}
	\subfigure{\includegraphics[height = 0.18\textwidth]{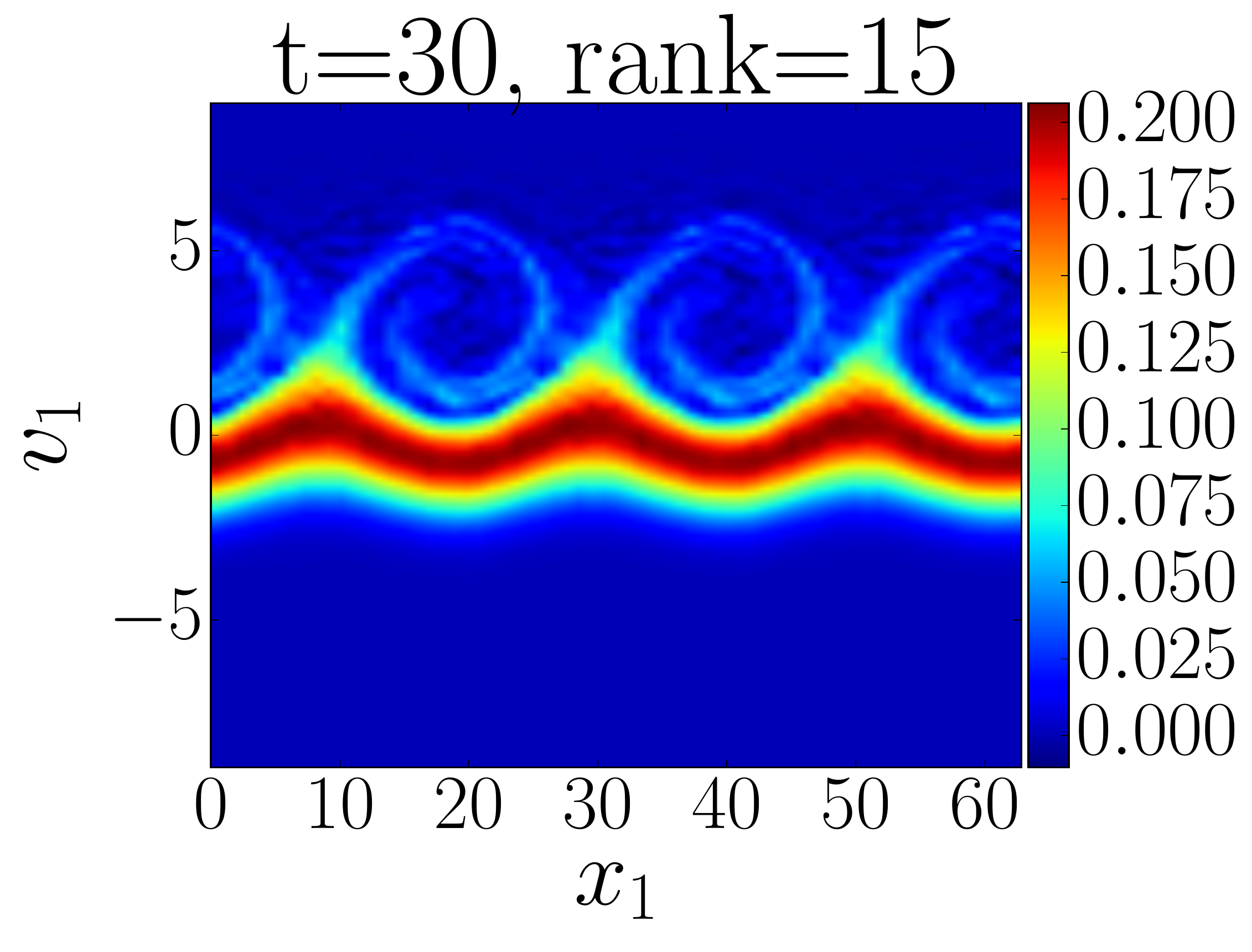}}
	\caption{Bump-on-tail instability, slices with $v_2=0$ of the distribution function $f$ at different times are shown for algorithm \ref{alg:Strang} (without correction) for $n_{x_1}=65$, $n_{v_1}=128$, $n_{v_2}=128$, and $\tau=0.1$. Upper row: results for rank = 5, lower row: results for rank = 15.}
	\label{fig:bot}
\end{figure}

The aim of this section is the study of the bump-on-tail instability. 
We choose an initial value similar to the one in \cite{Crouseilles2010DiscontinuousGS}:
\[
f(0,x_1,v_1,v_2) = \frac{1}{\sqrt{2}\pi} \left( \alpha \ee^{-v_1^2/2} + \beta \ee^{-2(v_1-4.5)^2} \left(1+\gamma \cos( k x_1)\right)  \right) \ee^{-v_2^2} .
\]
The following numerical results are obtained with $\alpha=9/10$, $\beta=2/10$, $\gamma = 0.03$ and $k= 0.3$. We consider the phase space $\Omega_{x_1} \times \Omega_v$ to be $[0,20\pi] \times [-9,9]^2$. The electric field is initialized according to Gauss' law and no initial magnetic field is considered. 
The perturbation of the equilibrium results in an instability which consists of three vortices traveling periodically in phase space. The numerical simulation of this phenomenon is challenging. The three vortices must be kept separate and the filamentation must be accurately resolved. We report our results for algorithm \ref{alg:Strang} (without correction) in figure \ref{fig:bot} where we compare the distribution function at different times. In the upper row, the results obtained for rank = 5 are shown. The qualitative behavior of the solution improves significantly when choosing rank = 15. In both configurations, the formation of the vortices is resolved very well.

\section{Numerical results for the Weibel instability}\label{sect:weibel}
In this section we consider an instance of the Weibel instability; see \cite{CEF15}.
The initial particle density is
\[ 
f(0,x_1,v_1,v_2) = \frac{1}{\pi v_{th}^2 \sqrt{T_r}} \ee^{-(v_1^2+v_2^2/T_r)/v_{th}^2}\left(1+\alpha \cos(kx_1)\right),
\]
where we choose $v_{th}=0.02$, $T_r = 12$, $k = 1.25$, and $\alpha=10^{-4}$. The magnetic field at $t=0$ is chosen as 
\[ 
B_3(0,x_1) = \beta \cos(kx_1)
\]
with $\beta = 10^{-4}$. The electric field $E_1$ is initialized by \eqref{eq:1+2Poisson}.
The space domain here is $\Omega_{x_1} = [0,2\pi/k]$, and the velocity domain is chosen as $\Omega_{v_1}\times \Omega_{v_2} = [-0.3,0.3]^2$. 

As explained in section \ref{sect:implementation}, we approximate the spatial derivatives in Fourier space. This results in a numerical scheme that does not have any dissipation. This is the correct behavior from a physical point of view. However, in some configurations this can cause numerical instabilities. In order to avoid this problem, we add some artificial dissipation into Maxwell's equations as follows
\begin{align} 
& \partial_t E_1(t,x_1) = -\int_{\Omega_{v}} v_1 f(t,x_1,v) \dd v + \varepsilon  \partial^2_{x_1} E_1(t,x_1) \label{eq:dampE1}\\
& \partial_t E_2(t,x_1) = -\partial_{x_1} B_3(t,x_1) - \int_{\Omega_{v}} v_2 f(t,x_1,v) \dd v + \varepsilon \partial^2_{x_1} E_2(t,x_1) \label{eq:dampE2}\\
& \partial_t B_3(t,x_1) = -\partial_{x_1} E_2(t,x_1) + \varepsilon \partial^2_{x_1} B_3(t,x_1). \label{eq:dampB3} 
\end{align}
By tuning the value of the parameter $\varepsilon$, we can adjust the strength of the numerical dissipation. We use $\varepsilon = 10^{-2} h^2$ in all simulations in this section, where $h$ is the spatial grid size. This is sufficient to remove the numerical instability. 
Such an approach is commonly used and dates back to \cite{vonNeumann}. The need of introducing artificial dissipation is due to the fact that here, contrary to the previous examples, the magnetic field is strong and thus non-physical oscillations are amplified by the numerical scheme. This is \textit{not} an issue with the low-rank approximation. Eulerian schemes based on spectral methods, e.g.~the scheme proposed in \cite{CEF15}, also require artificial dissipation. In the implementation, we treat the diffusion implicitly.

\begin{figure}
	\includegraphics[width = .5\textwidth]{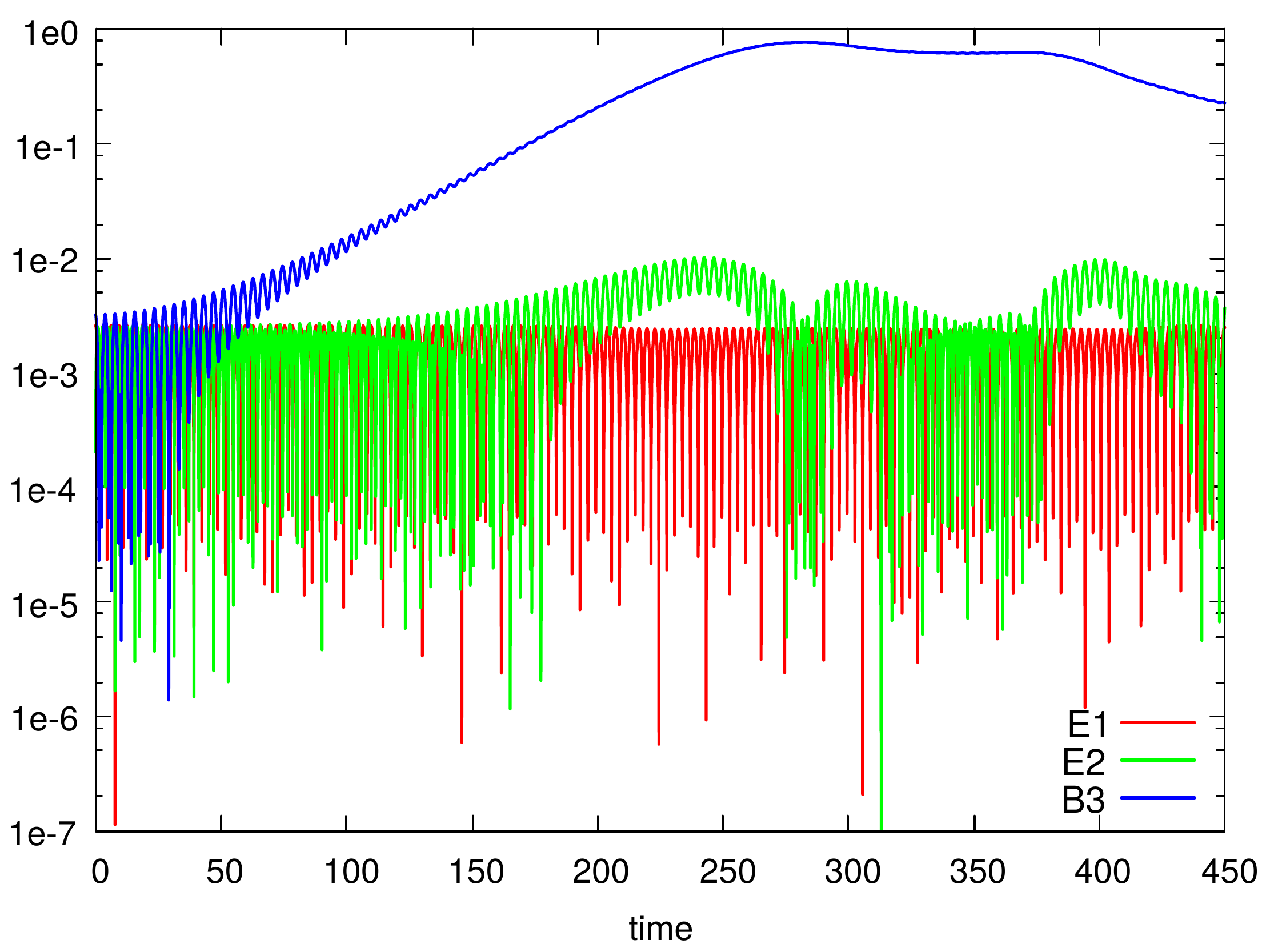}
	\includegraphics[width = .5\textwidth]{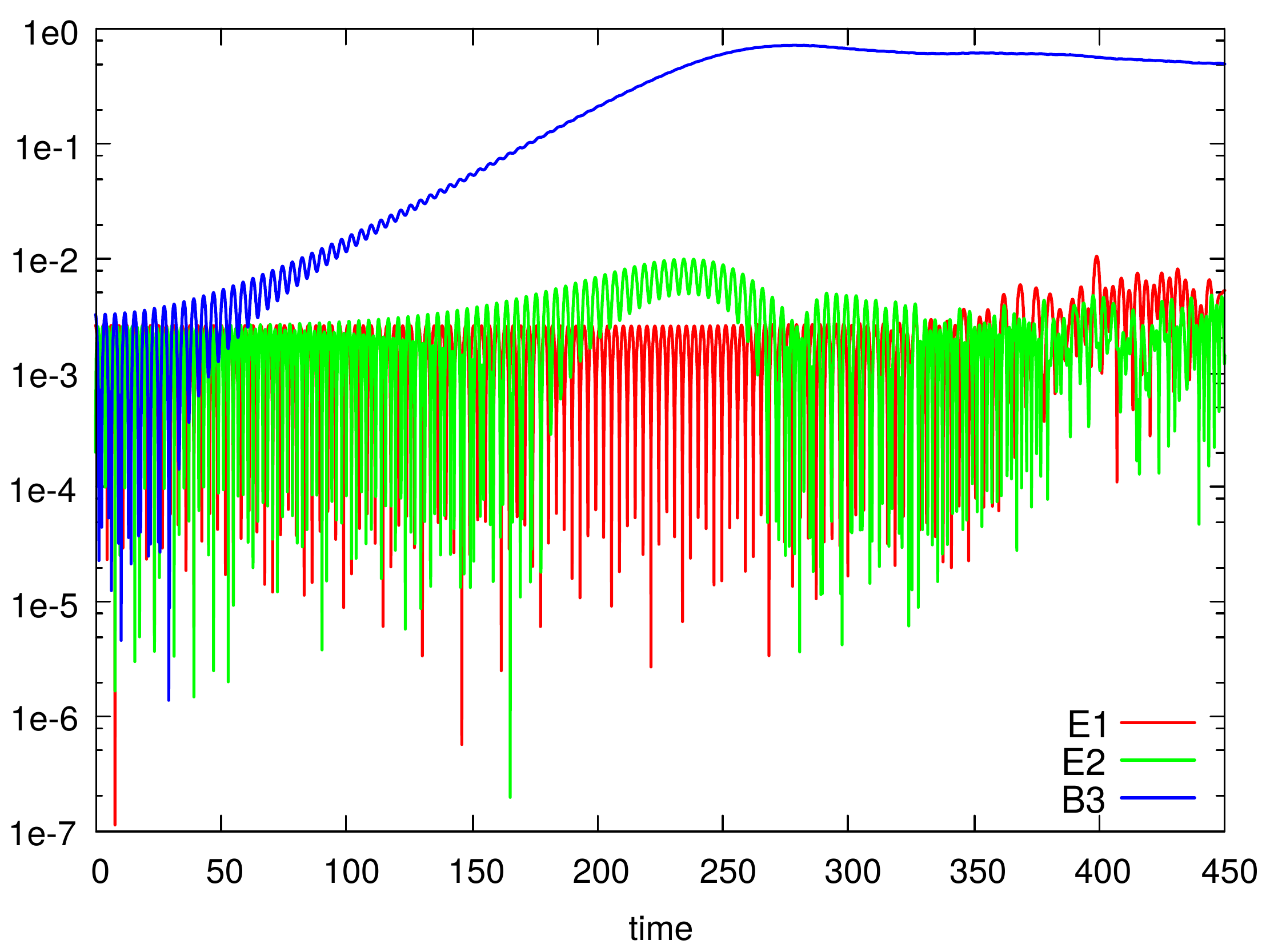}\\
	\center{
	\includegraphics[width = .5\textwidth]{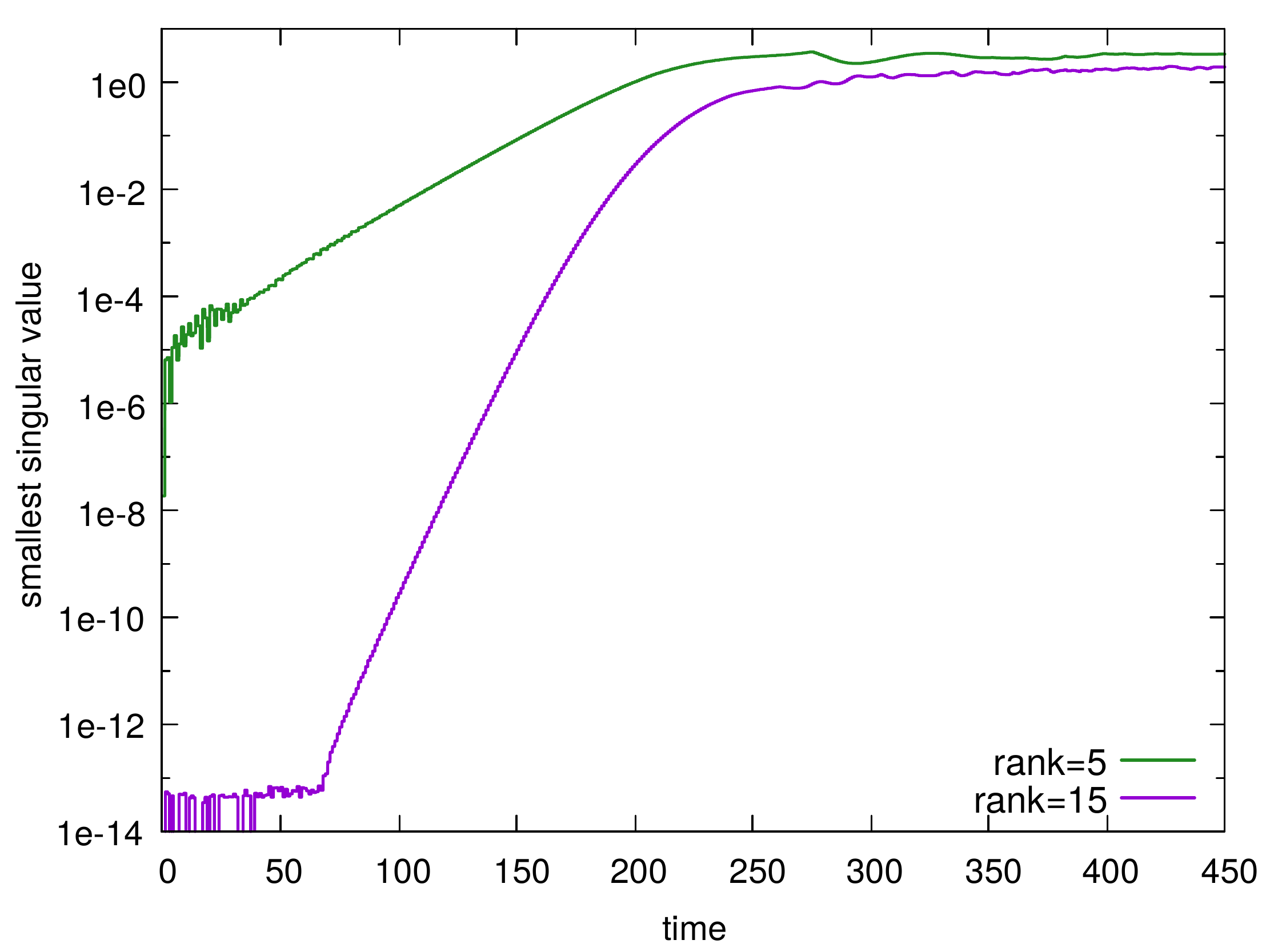}}
	\caption{Weibel instability. Results for algorithm \ref{alg:Strang} (without correction) for $n_{x_1}=65$, $n_{v_1}=192$, $n_{v_2}=192$, $\tau=0.05$, and $\varepsilon = 10^{-2}h^2$ (semi-log scale). The figures on the top show the first Fourier mode of each component of the electromagnetic field for rank = 5 and  rank = 15, respectively. The bottom figure displays the smallest singular value of the corresponding distribution function computed every 2 time steps.}
	\label{fig:wi_wc_energy}
\end{figure}

In figure \ref{fig:wi_wc_energy}, on the top, the first Fourier mode of the electric and magnetic fields, respectively, are plotted for algorithm \ref{alg:Strang} (without correction). We compare the results obtained for rank = 5 (left) and rank = 15 (right). We observe that the integrator works very well in both the linear regime and to simulate saturation. More specifically, it reproduces the exponential growth of the magnetic field until the saturation point. We refer to \cite[sect. 5.2]{CEF15} for some reference results and conclude that both choices of the approximation rank give good qualitative results also in the nonlinear regime.
On the other hand, the behavior of the smallest singular value (bottom right figure) indicates that the error in the distribution function is quite large in the nonlinear regime. Thus, we once again have a situation in which the averaged quantities are resolved much better than the distribution function itself.

\begin{figure}
	\includegraphics[width = .5\textwidth]{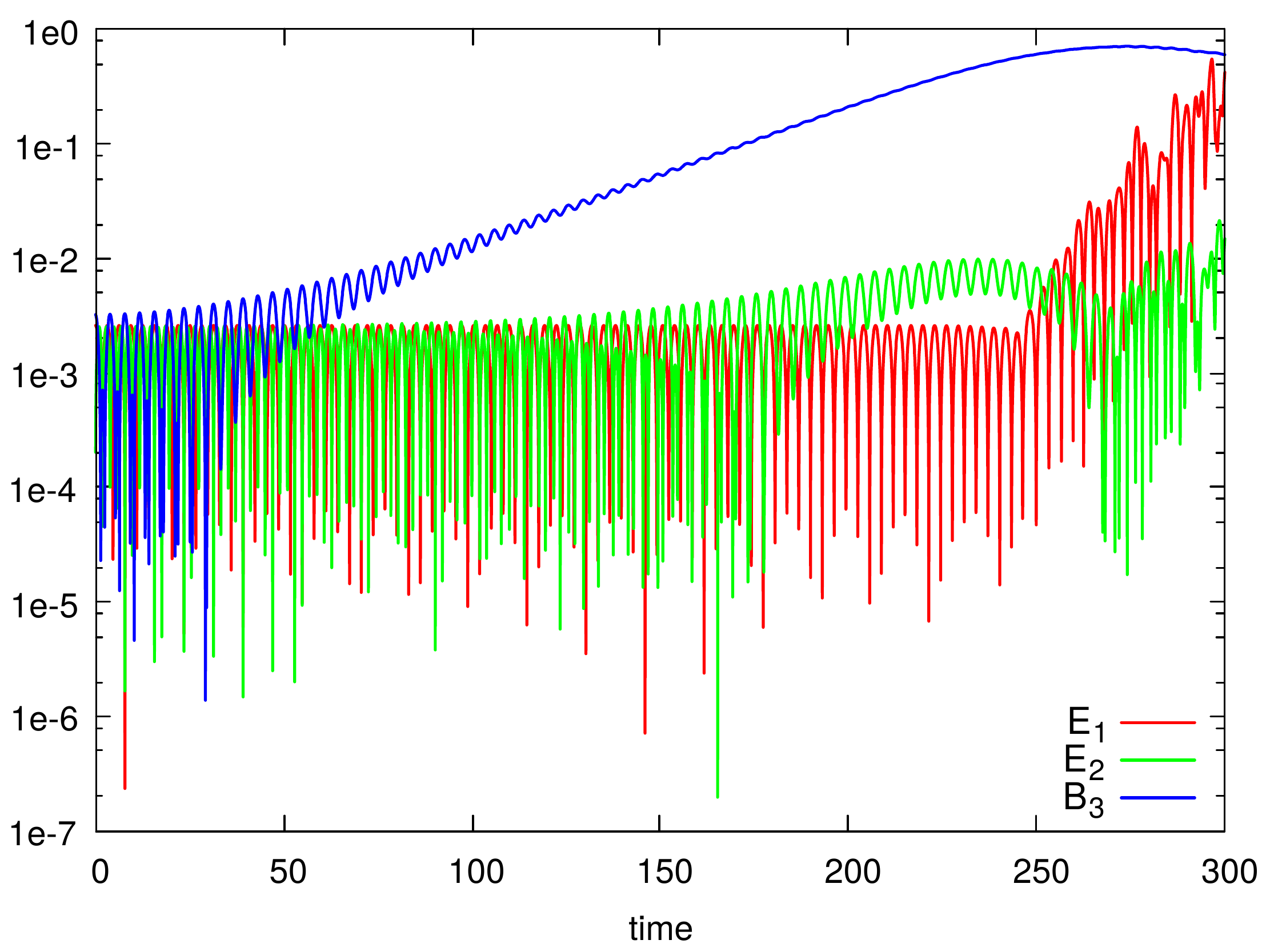} 
	\includegraphics[width = .5\textwidth]{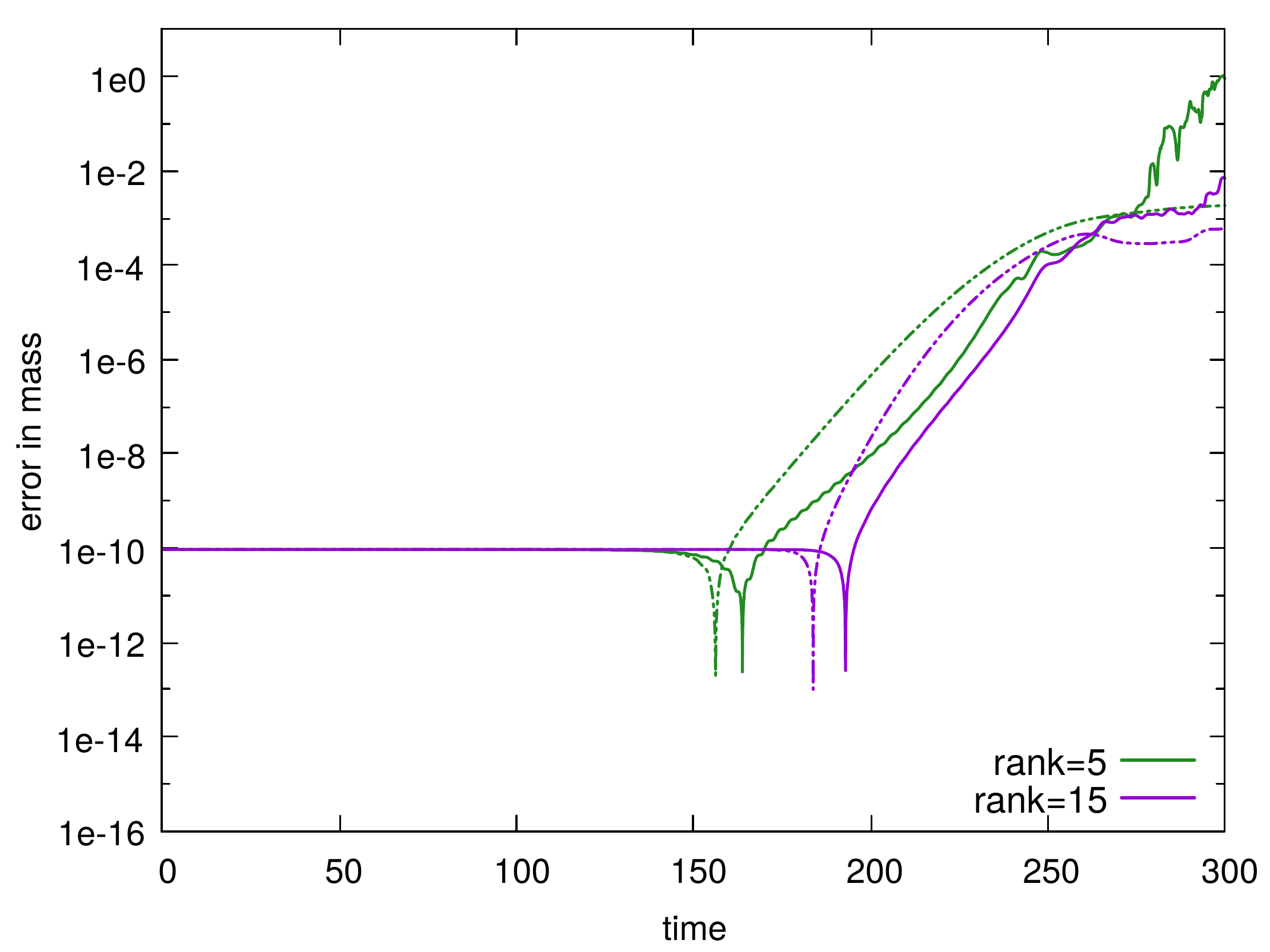} \\
	\includegraphics[width = .5\textwidth]{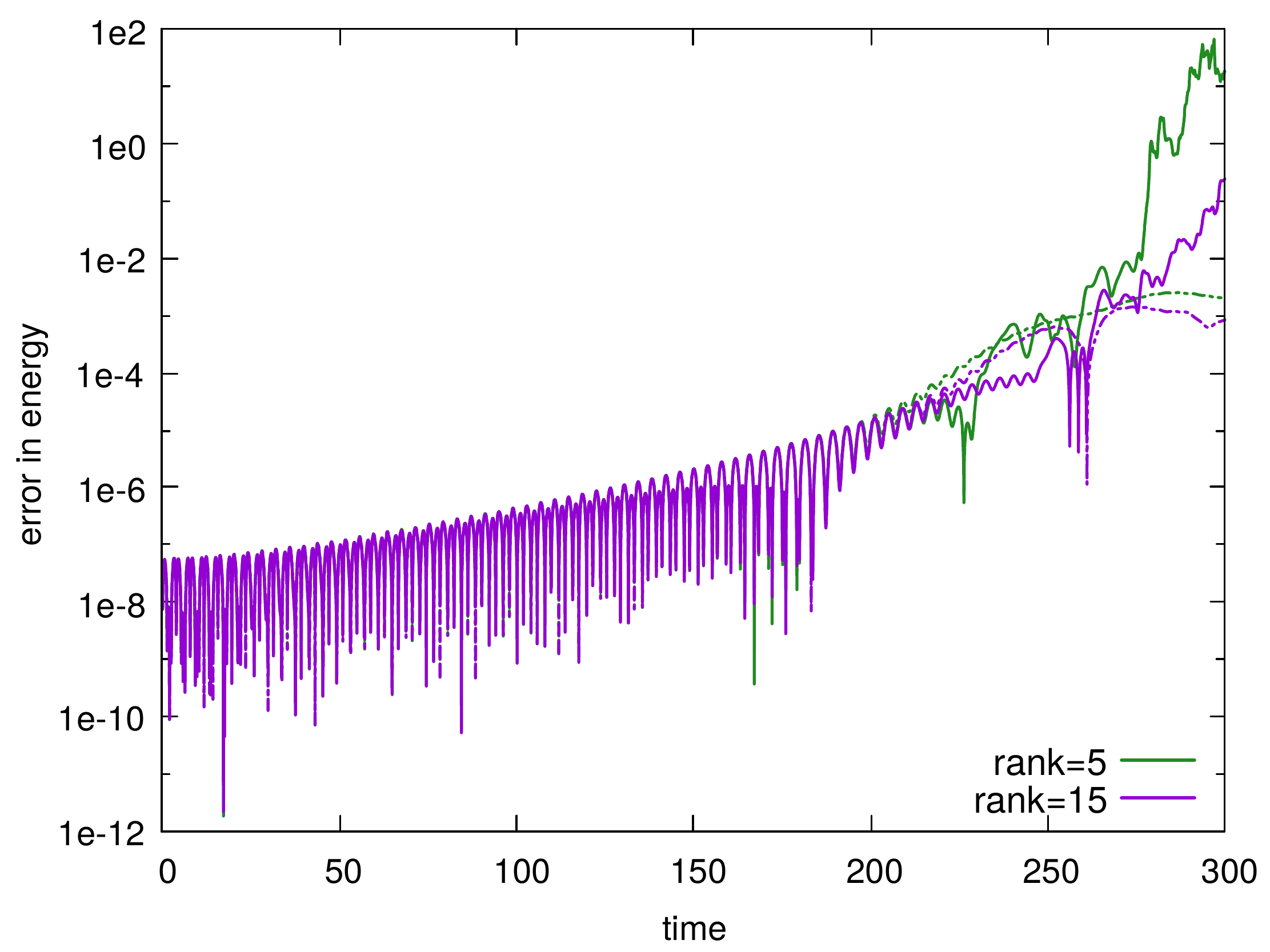} 
	\includegraphics[width = .5\textwidth]{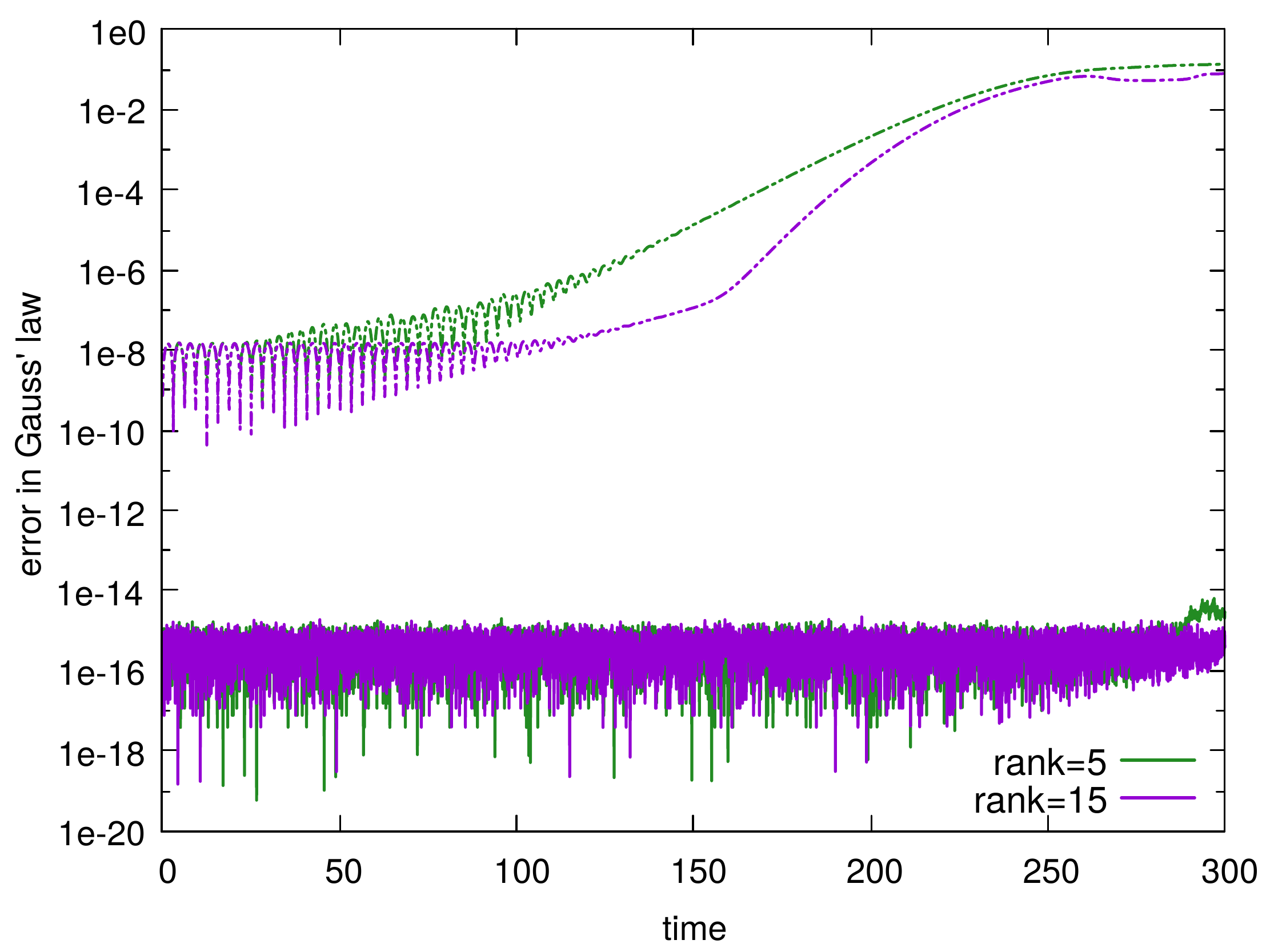}
	\caption{Weibel instability. Results for $n_{x_1}=65$, $n_{v_1}=192$, $n_{v_2}=192$, $\tau=0.05$, and $\varepsilon = 10^{-2}h^2$ (semi-log scale). The top left figure shows the first Fourier mode of each component of the electromagnetic field obtained with algorithm \ref{alg:Strang_corr} (with correction). In the other plots we compare the low-rank integrator described in algorithm \ref{alg:Strang} (without correction, dashed-dotted lines) with the divergence preserving scheme given in algorithm \ref{alg:Strang_corr} (with correction, full lines) in terms of conservation of mass (top right), energy (bottom left) and error in Gauss' law measured measured in the discrete $L^2$ norm (bottom right).}
	\label{fig:wi_comparisons}
\end{figure}

Similar conclusions can be drawn from figure \ref{fig:wi_comparisons}. Mass and energy are well conserved for both the schemes in the case of short time integration, but a significant error is incurred in the nonlinear regime. 
However, the divergence correction influences the overall behavior of the scheme in the nonlinear regime. In the top left figure we observe worse qualitative behavior of the Fourier modes, in particular for the electric field. Enforcing Gauss' law results in a deterioration of the quality of the other invariants we considered.
Let us explain this behavior. In general, there is no guarantee that a small correction is sufficient to obtain a point on the low-rank manifold which satisfies Gauss' law. Thus, there is a competition between the approximation error and the divergence correction. This is precisely what happens in the nonlinear regime. Therefore, introducing the divergence correction deteriorates the qualitative behavior of the numerical solution. This is primarily a problem in situations where the low-rank error in the distribution is large. 
We also remark that the divergence correction procedure explained in section \ref{sect:correction} has to be adapted to the damped equations for the electromagnetic field \eqref{eq:dampE1}--\eqref{eq:dampB3}.

\section*{Acknowledgements}
We would like to thank Bruno Despr\'{e}s (Laboratoire Jacques-Louis Lions) for the  helpful discussions.
We also thank the referees for their helpful comments, which
improved the presentation of this paper.

\bibliography{references}

\end{document}